\documentclass[11pt]{amsart}
\usepackage{palatino}
\usepackage[all]{xy}
\usepackage{amssymb}
\usepackage{latexsym}
\usepackage{amscd}
\usepackage{color}
\input amssym.def
\input amssym
\newtheorem{theorem}{Theorem}[section]
\newtheorem{lemma}[theorem]{Lemma}

\newtheorem{example}[theorem]{Example}
\newtheorem{proposition}[theorem]{Proposition}
\newtheorem{corollary}[theorem]{Corollary}
\newtheorem{conjecture}[theorem]{Conjecture}
\newtheorem{remark}{Remark}

\newcommand{\trip}{\mathcal{W}(p)}
\newcommand{\sing}{\overline{M(1)_p}}
\newcommand{\no}{\noindent}
\def \<{\langle}
\def \>{\rangle}

\def \a{\alpha }

\def \l{\lambda }
\def \ga{\gamma }

\newcommand{\bea}{\begin{eqnarray}}
\newcommand{\eea}{\end{eqnarray}}
\newcommand{\be}{\begin {equation}}
\newcommand{\ee}{\end{equation}}

\newcommand{\wt}{{\rm {wt} }   }

\newcommand{\Z}{\Bbb Z}

\newcommand{\W}{\mathcal W}
\newcommand{\mip}{\overline{M(1)_p}}
\newcommand{\Zp}{{\Bbb Z}_{>0} }

\newcommand{\N}{{\Bbb Z}_{\ge 0} }
\newcommand{\C}{\Bbb C}

\newcommand{\la}{\langle}
\newcommand{\ra}{\rangle}

\newcommand{\halmos}{\rule{1ex}{1.4ex}}

\newcommand{\epfv}{\hspace*{\fill}\mbox{$\halmos$}}

\usepackage{amssymb}

\begin{document}
\title[Triplet vertex algebra
$\mathcal{W}(p)$]{On the triplet vertex algebra $\mathcal{W}(p)$}
\author{Dra\v{z}en Adamovi\'c and Antun Milas}

\begin{abstract}
\noindent We study the triplet vertex operator algebra
$\mathcal{W}(p)$ of central charge $1-\frac{6(p-1)^2}{p}$, $p \geq
2$. We show that $\trip$ is $C_2$-cofinite but irrational since it
admits indecomposable and logarithmic modules. Furthermore, we
prove that $\trip$ is of finite-representation type and we provide
an explicit construction and classification of all irreducible
$\mathcal{W}(p)$-modules and describe block decomposition of the
category of ordinary $\trip$-modules. All this is done through an
extensive use of Zhu's associative algebra together with explicit
methods based on vertex operators and the theory of automorphic
forms. Moreover, we obtain an upper bound for ${\rm
dim}(A(\mathcal{W}(p)))$. Finally, for $p$ prime, we completely
describe the structure of $A(\trip)$. The methods of this paper
are easily extendable to other $\mathcal{W}$-algebras and
superalgebras.

\end{abstract}


\address{Department of Mathematics, University of Zagreb, Croatia}
\email{adamovic@math.hr}

\address{Department of Mathematics and Statistics,
University at Albany (SUNY), Albany, NY 12222}
\email{amilas@math.albany.edu}

\maketitle

\renewcommand{\theequation}{\thesection.\arabic{equation}}
\renewcommand{\thetheorem}{\thesection.\arabic{theorem}}
\setcounter{equation}{0} \setcounter{theorem}{0}
\setcounter{equation}{0} \setcounter{theorem}{0}
\setcounter{section}{-1}

\section{Introduction}

The main focus of vertex operator algebra theory so far has been
on understanding {\em rational vertex operator algebras}.
%\footnote{A
%vertex operator algebra $V$ is usually called rational if it has
%finitely many equivalence classes of irreducible modules, such that
%every $V$-module is completely reducible.}.
This progress has led,
in particular, to several important breakthroughs in the area such
as Zhu's modular invariance theorem \cite{Zhu} and Huang's recent
proof of the Verlinde's conjecture \cite{H1}. All these developments
are deeply rooted in ideas of conformal field theory.

On the other hand, irrational vertex operator algebras did not
attract so much attention for several obvious reasons. The
category of modules of a (irrational) vertex operator algebra has
often too many irreducible objects, which forces many important
features such as modular invariance to be absent. In view of that,
it is reasonable to focus first on vertex algebras of {\em
finite-representation type} and relax the semisimplicity
condition. A vertex operator algebra $V$ will be of
finite-representation type if its Zhu associative algebra is
finite-dimensional, which is the case when $V$ is $C_2$-cofinite
(cf. \cite{DLM}). Since the $C_2$-cofiniteness property is rather
strong, it is interesting that there are examples of irrational
$C_2$-cofinite vertex operator algebras (cf. \cite{Abe}). More
surprisingly, there is a version of modular invariance theorem for
modules of $C_2$-cofinite vertex algebras \cite{Miy}, so
irrational $C_2$-cofinite vertex algebras appear to be very
special objects. Since there is only a handful of known examples
of vertex algebras with these properties, it is important to study
the known examples and to seek for new models.

Our present line of work, continuing \cite{AdM}, is concerned with
"$\mathcal{W}$-algebras". The $\mathcal{W}$-algebras have been
studied intensively by the mathematicians and physicists over the
last two decades. These "algebras" are not Lie algebras in the
classical sense, but rather close cousins of vertex algebras
associated to affine Lie algebras and lattice vertex algebras, and
are usually defined via the process of {\em quantum reduction}
\cite{FKRW}, \cite{FrB}. The $\mathcal{W}$-algebras and their
classical counterparts are important not only {\em per se}, but
also in connection with integrable hierarchies, opers and even
geometric Langlands correspondence \cite{FrB}. For some recent
advances in the theory of $\mathcal{W}$-algebras see \cite{Ar},
\cite{dSK}, \cite{Fr}, etc.

The $\mathcal{W}$-algebras associated to affine Lie algebras come
in families parameterized by the value of central charge. For
generic values of the parameter, $\mathcal{W}$-algebras have
fairly explicit description and are known to be irrational (cf.
\cite{Ar}, \cite{FKRW}). For non-generic values the situation
seems to be opposite. The simplest examples of non-generic
$\mathcal{W}$-algebras - the vacuum Virasoro minimal models - are
known to be rational. But in higher "rank", there are almost no
classification results of representations of non-generic
$\mathcal{W}$-algebras. A sole example would be the proof of
rationality of a rank two $\mathcal{W}_3$-algebra of central
charge $c=\frac{6}{5}$, which is already quite involved
\cite{DLTYY}. Even so, one does expect that, suitably defined,
vacuum minimal models for higher rank $\mathcal{W}$-algebras are
also rational vertex algebras \cite{FrB}.

In this paper we focus on a prominent family of
$\mathcal{W}$-algebras - called {\em triplet}
$\mathcal{W}$-algebras - introduced by Kausch in \cite{Ka1},
\cite{Ka2}. These are parameterized by a positive integer $p \geq
2$, and have central charge $c=1-\frac{6(p-1)^2}{p}$, $p \geq 2$.
Unlike the $\mathcal{W}$-algebras discussed in the previous
paragraph, the triplet is associated to a (non-root) lattice
vertex algebra. In terms of generators the triplet comes with four
distinguished vectors; one generator is the usual Virasoro vector
and the remaining three are certain primary fields of conformal
weight $2p-1$. For $p=2$, the triplet has central charge $-2$ and
it admits a realization via the so-called {\em symplectic
fermions} \cite{Ka1}, \cite{GK1}, \cite{GK2}. A great deal of
research has been done on the symplectic fermions from several
different points of view. One particularly interesting feature of
this model is the appearance of the so-called {\em logarithmic
modules} (i.e., modules with a non-diagonalizable action of the
Virasoro generator $L(0)$). This links the $c=-2$ model and more
general triplets with {\em logarithmic conformal field theory}
(LCFT), a new physical theory with possible applications in
condensed matter physics and string theory (cf. see
\cite{HLZ1}-\cite{HLZ2} for a vertex operator algebra approach to
LCFT). There is a large body of work devoted to various
interaction between the triplet $\mathcal{W}$-algebras, LCFT and
the quantum group $\mathcal{U}_q(\goth{sl}_2)$ at a root of unity.
For these and other related developments we refer the reader to
\cite{EFH}, \cite{FFHST}, \cite{FHST}, \cite{FGST1}-\cite{FGST3},
\cite{F1}, \cite{FG}, \cite{GK1}, \cite{GK2}, \cite{GR},
\cite{Ka1}, \cite{Ka2}, \cite{KaW}, and especially two excellent
reviews: \cite{F3} and \cite{Ga}.

Since the triplet $\mathcal{W}$-algebras are in fact vertex
(operator) algebras, it is tempting to use vertex algebra theory to
analyze these objects. In \cite{Abe} Abe studied the $c=-2$ triplet
vertex algebra, denoted by $\mathcal{W}(2)$, by using vertex algebra
theory. He eventually succeeded in classifying irreducible
$\mathcal{W}(2)$-modules and describing the Zhu's associative
algebra $A(\mathcal{W}(2))$. A rather explicit fermionic
construction of $\mathcal{W}(2)$ played a prominent rule in Abe's
approach. However, for $p>2$ the triplet vertex algebra $\trip$ has
no natural fermionic realization so one needs a completely different
approach (cf. \cite{FHST}) to study these vertex algebras and to
extend classification results to $\trip$, for $p > 2$. That was
precisely the main motivational problem for this paper.

Here is a short description of our main results. We start with a
rank one lattice vertex algebra $V_L$ and view the triplet $\trip$
as the kernel of a screening operator acting on $V_L$. It is
important to say here that $L$ is not a root lattice.
Alternatively, the triplet can be defined in terms of generators
as we indicated earlier. Our first two results are the
classification of irreducible $\trip$-modules and a fairly simple
proof of the $C_2$-cofiniteness:

\newpage

\noindent {\bf Theorem A.} {\em \begin{itemize} \item[(i)] The
triplet vertex algebra $\trip$ is $C_2$-cofinite \item[(ii)] The
triplet vertex algebra $\trip$ has precisely $2p$ inequivalent
irreducible modules. \item[(iii)] The vertex algebra $\trip$ is
irrational.
\end{itemize}
}
\vskip 2mm

For precise statements see Theorem \ref{c2}, Theorem
\ref{ireducibilni-moduli} and Theorem \ref{class-ired-modules}.
One should say that Theorem A is in agreement with the results
obtained by physicists. We also stress that the $C_2$-cofiniteness
part (i) was discussed in \cite{CF} by using a completely
different approach.

Since the triplet vertex algebra is not rational (cf. Proposition
\ref{ses}) it is not evident how to obtain a decomposition of
Zhu's associative algebra $A(\trip)$ as a direct sum of its
ideals. Even though Theorem A indicates that $A(\trip)$ should be
a sum of $2p$ ideals, it is not clear what precisely are those
ideals and how to compute their dimensions. Our next result gives
a partial answer to this problem.

\vskip 2mm

\noindent {\bf Theorem B.} {\em  Zhu's algebra $A(\trip)$
decomposes as a direct sum of ideals
$$A(\trip)=\bigoplus_{i=2p}^{3p-1}  \mathbb{M}_{h_{i,1}} \oplus \bigoplus_{i=1}^{p-1}
\mathbb{I}_{h_{i,1}} \oplus \mathbb{C},$$ where
$\mathbb{M}_{h_{i,1}}$ is an ideal isomorphic to
$M_2(\mathbb{C})$, and each $\mathbb{I}_{h_{i,1}}$ is at most
two-dimensional (all ideals are parameterized by certain conformal
weights $h_{i,1}$).}

\vskip 2mm

In fact, we prove much more; all our ideals are described with
explicit spanning sets (bases?), which is useful for computational
purposes.

Notice that the last result is not completely satisfactory because
we really do not know whether $\mathbb{I}_{h_{i,1}}$ is
two-dimensional (this is our conjecture though). The problem
relies on the existence of certain logarithmic $\trip$-modules,
predicted by physicists (cf. \cite{FHST}). In the $p=2$ case the
existence of such a module can be easily seen be using symplectic
fermions, but for $p>2$ it is not quite clear how to construct
these modules explicitly. Nevertheless, by using a work of
Miyamoto \cite{Miy} we settle this problem, at least when $p$ is a
prime integer.

\vskip 2mm

\noindent {\bf Theorem C.} {\em Suppose that $p$ is a prime
number. Then each $\mathbb{I}_{h_{i,1}}$ is two-dimensional and
$${\rm dim}(A(\trip))=6p-1.$$ }

\vskip 2mm

For general $p$ we only have a partial result in this direction
(cf. Proposition \ref{gen-log-p}).

Besides an obvious representation theoretic interest, our results
and techniques have other merits. For example, methods used in
this paper are applicable to other $\mathcal{W}$-algebras and
superalgebras defined via screening operators. This paper also
gives rigorous proofs of several claims in physics literature
about the triplet, their representations and logarithmic modules.
Thus, our results and have immediate applications in logarithmic
conformal field theory.

\newpage

\section{The triplet vertex algebra $\mathcal{W}(p)$}

In this section we introduce the triplet vertex algebra and study
some of its representations.

The setup is similar as in \cite{AdM} (see also \cite{A-2003}) so
we omit many details. For an integer $p \geq 2$, we fix a rank one
lattice $\mathbb{Z}\alpha$, generated by $\alpha$ with
$$\langle \alpha, \alpha \rangle=2p.$$ We denote by $(V_L,Y,\omega,{\bf 1})$ the
corresponding lattice vertex operator algebra \cite{D}, \cite{LL}.
As a vector space, $$V_{L}=\mathcal{U}(\hat{\goth{h}}_{<0})
\otimes
 \mathbb{C}[L],$$
where $\mathbb{C}[L]$ is the group algebra of $L$ and
$\hat{\goth{h}}$ is the affinization of the one-dimensional
algebra spanned by $\alpha$, the vacuum vector
$${\bf 1}=1 \otimes 1$$
and the conformal vector
$$\omega=\frac{\alpha(-1)^2}{4p}{\bf 1}+\frac{p-1}{2p}\alpha(-2){\bf 1}.$$
Notice that this is not the usual quadratic Virasoro generator
used throughout the literature.
%
% isp. 4p u 2p
%
For convenience, let us fix $h=\frac{\alpha}{\sqrt{2p}}$, so that
$\langle h, h \rangle=1$.

It is known that the vertex operator algebra $V_L$ is rational, in
the sense that it has finitely many irreducible $V_L$-modules and
that every $V_L$-module is completely reducible (see \cite{LL} for
instance). If we denote by $L^\circ=\frac{\mathbb{Z} \alpha}{2p}$
the dual lattice of $L$, then
$$V_{L+\lambda}=\mathcal{U}(\hat{\goth{h}}_{<0}) \otimes
e^{\lambda} \mathbb{C}[L], \ \ \ \lambda \in L^\circ/L,$$ are, up
to equivalence, all irreducible $V_L$-modules. For $\lambda=0$ we
recover the vertex algebra $V_L$, while $V_{L^\circ}$ has a
generalized vertex operator algebra structure. The Virasoro
algebra acts on $V_{L+\lambda}$ with the central charge
$$c_{p,1}=1-\frac{6(p-1)^2}{p}.$$
As usual the Virasoro generators will be denoted by $L(n)$, $n \in
\mathbb{Z}$. Then the degree zero operators $L(0)$ and the charge
operator $h(0)$ equip $V_{L+\lambda}$ with a compatible bigrading.

We also define rational numbers
$$h_{m,n}=\frac{(m-np-p+1)(m-np+p-1)}{4p}, \ \ m,n \in \mathbb{Z}$$
parameterizing $(1,p)$-minimal models models at the boundary of
Kac's table. In fact, it is sufficient to consider the weights \be
\label{cp1wt} h_{m+1,1}=\frac{m(m-2p+2)}{4p}, \ \ m \geq 0. \ee

As in \cite{A-2003} and  \cite{AdM} we let \bea && Q : V_L
\longrightarrow V_L, \nonumber \\
&& \tilde{Q} : V_L \longrightarrow V_{L^\circ},\nonumber \eea
where \bea && Q=e^{\alpha}_0, \nonumber \\ &&
\tilde{Q}=e^{-\alpha/p}_0 \nonumber \eea are {\em screening}
operators introduced by Dotsenko and Fateev (cf. \cite{FrB}). As
usual, here
$$Y(e^{\alpha},x)=\sum_{n \in \mathbb{Z}} e^{\alpha}_{n} x^{-n-1}, \ \ \ Y(e^{-\alpha/p},x)=\sum_{r \in \tfrac{1}{p}\mathbb{Z}}
e^{-\alpha/p}_{r} x^{-r-1},$$ where the second vertex operator
belong to the generalized vertex algebra $V_{L^\circ}$. We stress
here that the screening operator $Q$ acts as a derivation
$$Q (a_n b)=(Qa)_n \ b + a_n \ (Qb), \ \ a,b \in V_L, n \in \mathbb{Z}.$$
This formula will be extensively used throughout the paper.

Inside the vertex algebra $V_L$ we consider the following three
vectors and their vertex operators, crucial for our present work:
\bea F&=&e^{-\alpha}, \nonumber \\
H&=&Qe^{-\alpha},\nonumber \\
E&=&Q^2 e^{-\alpha}, \nonumber \eea
$$Y(X,z)=\sum_{i \in \mathbb{Z}} X_{i} z^{-i-1}, \ \ X \in \{E,F,H
\}.$$ As we shall see later, the choice of letters $F$, $H$ and
$E$ resembling the standard $\goth{sl}_2$ generators is not an
accident.

Let us recall the formulas (cf. \cite{A-2003}, \cite{AdM}): \bea
&&\label{EFH} E_i E=F_i F
=Q(H_i H)=0, \ \ i \geq -2p. \\
\label{comut} && [Q, \tilde{Q}] = 0.
\eea

We shall need the following useful formulae which hold in $V_L$

 \bea Y(e^{\a},x_1)  Y(e^{\a},x_{2}) & =&
 E^{-}(-\a, x_1,  x_{2}) E^{+}(-\a, x_1,  x_{2})
 (x_1 -x_{2}) ^{2p} e ^{2 \a}(x_1 x_{2} ) ^{\a}
 \label{prod-vo}
\eea
where
$$ E^{\pm}(-\a, x_1,  x_{2}) = \mbox{exp} \left( \sum_{k=1}
^{\infty}\tfrac{\alpha(\pm k)}{\pm k}( x_1 ^{\mp k}  + x_{2} ^{\mp
k})\right).$$

For $i \in {\Z}$, we  set
\bea && \ga_i = \frac{i}{2p} \a \label{ga-def}. \eea

We shall first present results on the structure of $V_L$--modules
as modules for the Virasoro algebra. By using Lemma 4.3 from
\cite{AdM}  and  the structure theory of Feigin-Fuchs modules
\cite{FF}  we get the following theorems. (Here we also use that
every $V_L$--module is a direct sum of Feigin-Fuchs modules due to
charge decomposition.)

%and
%a completely analogous proof to those of Theorem \ref{simple} and Theorem \ref{irr-module}
%we get the following result.

\begin{theorem} \label{str-ff-1}
Assume that $i \in \{0, \dots, p-2\}$.
 \item[(i)] As a Virasoro algebra module, $V_{L+\ga_i}$  is generated by the family of
singular and cosingular vectors $ \widetilde{Sing}_{i}  \bigcup
\widetilde{CSing}_{i}$, where
$$  \widetilde{Sing}_{i} =  \{ u_i ^{(j, n)} \ \vert \ j, n \in {\N}, \ 0 \le j \le 2n \}; \
\
  \widetilde{CSing}_{i} =  \{ w_i ^ {(j, n)} \ \vert \ n \in {\Zp}, 0 \le j \le 2n +1 \}. $$
These vectors satisfy the following relations:
\bea
&&  u_i ^{(j, n)}= Q ^{j} e ^{\ga_i -n \alpha} , \ \ Q ^{j} w_i
^{(j, n)} = e ^{\ga_i + n \a}. \nonumber \eea
 The submodule generated by singular vectors
$\widetilde{Sing}_i$
 is isomorphic to $$
%\Lambda(i+1)
\overline{V_{L+\ga_i}} \cong \bigoplus_{n =0 } ^{\infty} (2n +1)
L(c_{p,1}, h_{i+1, 2n+1}).$$

\item[(ii)]The quotient module  is isomorphic to
$$V_{L+\ga_i} / \overline{V_{L+\ga_i}}   \cong
%\Pi( p-i-1) \cong
  \bigoplus_{n =1 } ^{\infty}
(2 n ) L(c_{p,1}, h_{i+1, -2n +1}). $$

\item[(iii)] As a Virasoro algebra module $V_{L + \ga_{p-1}}$ is
generated  by the family of singular vectors
$$  \widetilde{Sing}_{p-1} =  \{ u_{p-1} ^{(j, n)} : = Q ^{j} e ^{\ga_{p-1} -n \alpha}  \ \vert \ j, n \in {\N}, \ 0 \le j \le 2n \}; $$

 and it is isomorphic to
$$
%\Lambda(p)
 V_{L + \ga_{p-1}} \cong  \bigoplus_{n=0} ^{\infty} (2n +1)L(c_{p,1}, h_{p, 2n
+1}).$$

\end{theorem}

\begin{theorem} \label{str-ff-2}
Assume that $i \in \{p, \dots, 2 p-2\}$.
 \item[(i)] As a Virasoro algebra module, $V_{L+\ga_{i}}$  is generated by the family of
singular and cosingular vectors $ \widetilde{Sing}_{i}  \bigcup
\widetilde{CSing}_{i}$, where
$$  \widetilde{Sing}_{i} =  \{ u_i ^{(j, n)} \ \vert \ j\in {\N},  n \in {\Zp}, \ 0 \le j \le 2n-1 \}; \
\
  \widetilde{CSing}_{i} =  \{ w_i ^ {(j, n)} \ \vert \ j, n \in {\N}, 0 \le j \le 2n \}. $$
These vectors satisfy the following relations:
\bea
&&  u_i ^{(j, n)}= Q ^{j} e ^{\ga_i -n \alpha} , \ \ Q ^{j} w_i
^{(j, n)} = e ^{\ga_i + n \a}. \nonumber \eea
 The submodule generated by singular vectors
$\widetilde{Sing}_i$
 is isomorphic to $$
%\Pi( i+1 - p)
\overline{V_{L+\ga_i}} \cong \bigoplus_{n =1 } ^{\infty} (2n)
L(c_{p,1}, h_{i+1, 2n+1}).$$

\item[(ii)]The quotient module  is isomorphic to
$$V_{L+\ga_i} / \overline{V_{L+\ga_i}}  \cong
% \Lambda( 2 p-i-1) \cong
\bigoplus_{n =0} ^{\infty} (2 n +1) L(c_{p,1}, h_{i+1, -2n +1}).
$$

\item[(iii)] As a Virasoro algebra module $V_{L + \ga_{2 p -1}}$
is generated  by the family of singular vectors
$$  \widetilde{Sing}_{2 p-1} =  \{ u_{2 p-1} ^{(j, n)} : = Q ^{j} e ^{\ga_{2p-1} -n \alpha}  \ \vert \  n \in {\Zp}, j \in {\N},  \ 0 \le j \le 2n-1 \}; $$

 and it is isomorphic to
$$
%\Pi(p)
V_{L+\ga_{2p-1}} \cong  \bigoplus_{n=1} ^{\infty} (2n)L(c_{p,1},
h_{2 p, 2n +1}).$$

\end{theorem}

We will be concerned with certain (vertex) subalgebras of the
lattice vertex algebra $V_L$. For these purposes we recall a few
basic definitions. Let $V$ be a vertex (operator) algebra and $S
\subset V$. Then we denote by $\langle S \rangle$ the vertex
(operator) subalgebra generated by $S$ (i.e., the smallest vertex
(operator) algebra containing the set $S$). Similarly if $W$ is a
subset of a $V$-module $M$, we denote by $\langle W \rangle$ the
submodule of $M$ generated by $W$.

The vertex (operator) algebra $\langle S \rangle$ is said to be
{\em strongly generated} (cf. \cite{Kac}) by $S$ if it is spanned
by vectors of the form
$$(v_1)_{n_1}(v_2)_{n_2} \cdots (v_k)_{n_k}{\bf 1}, \ \ v_i
\in S, \ \ n_i <0.$$

Now, we are ready to introduce the main protagonist of the paper.
For $p$ as above we denote by $\mathcal{W}(p) \subset V_L$ the
following vertex operator algebra (cf. \cite{FHST})
$$\mathcal{W}(p):={\rm Ker}|_{V_L} \tilde{Q}.$$

Notice that this resemble the definition of the singlet vertex
algebra $\sing$ studied in \cite{A-2003} and \cite{AdM}, with the
only difference that now the kernel is taken over the whole
lattice algebra, whereas the singlet vertex algebra is the kernel
of $\tilde{Q}$ restricted onto the charge zero subalgebra $M(1)
\subset V_L$. We call
$\trip$ the {\em triplet vertex algebra}. Clearly, $\mip$ is a
vertex subalgebra of $\trip$. Recall also that $\mip$ is a simple
vertex opeator algebra strongly generated by $\omega$ and $H$ (cf.
\cite{A-2003}, \cite{AdM}). Next result will identify the
generators for $\trip$.

\begin{proposition} \label{stgen}
\item[(i)] We have $\trip = \overline{V_L}$. \item[(ii)] The
vertex operator (sub)algebra $\trip \subset V_L$ is strongly
generated by $E,H,F$ and $\omega$.
\end{proposition}
\noindent {\em Proof.}
Recall the structure of $V_L$ as a module for the Virasoro algebra
from Theorem \ref{str-ff-1}. By using (\ref{comut}), similarly to
the proof of Theorem 3.1 in  \cite{A-2003},  we conclude  that
$\trip$ is a completely reducible  module for the Virasoro algebra
which is generated  by the family of singular vectors:
\bea \label{sing-trip} && Q^{j} e^{-n\a}, \quad n \in {\N}, \ j
\in \{0, \cdots, 2n \}. \eea
This proves (i).  Let $W_n$ be the Virasoro module generated by
singular vectors
$$Q^{j} e^{-m\a}, \quad m \le n, \ j \in {\N} .$$
Therefore $\trip = \bigcup_{n\in {\N} } W_n$.
Let now $U$ be the subspace  of $\trip$  spanned by vectors of the
form
$$(v_1)_{n_1}(v_2)_{n_2} \cdots (v_k)_{n_k}{\bf 1}, \ \ v_i
\in \{E,F,H,\omega\} \ \ n_i <0,$$
 Clearly, $U \subseteq \trip$.
We shall prove that $U=\trip$, i.e., $\trip$ is strongly generated
by $E,F,H, \omega$.

 In order to
prove that $U = \trip$ it is enough to show that  $W_n \subseteq
U$ for every $n \in {\Zp}$. We shall prove this claim by induction
on $n$. By the definition, the claim holds for $n=1$. Assume now
that $W_n \subseteq U.$ Set $j_0 = 2n p +1$. The results from
\cite{A-2003} and \cite{AdM} imply that
\bea
 F_{-j_0}  e^{-n \a} &=& e^{-(n+1) \a}, \nonumber \\
 E_{-j_0} Q^{2n}
e^{-n\a} &=& C_{2n+1} Q^{2n +2} e^{-(n+1)\a}, \nonumber \eea {where}
$C_{2n+1} \ne 0$ and
$$  H_{-j_0} Q^{j} e^{-n\a} = C_j  Q^{j+1} e^{-(n+1)\a} + v_j ' ,
$$
where $ v_j' \in W_n, \ C_j \neq 0, \ 1 \le j \le 2n.$
These
relations imply that $W_{n+1} \subseteq U$. By induction we conclude
that $W_n \subseteq U$ for every $n \in {\Zp}$ and therefore $U =
{\trip}$. \epfv

\vskip 5mm

It is not hard to see that in fact
\begin{corollary}
The triplet vertex algebra $\trip$ is spanned by
$$L(-r_1) \cdots L(-r_k)H_{-s_1} \cdots H_{-s_l}X_{-t_1} \cdots
X_{-t_m}{\bf 1},$$ where $X=E$ or $X=F$, $k,l,m \in \mathbb{N}$,
$s_i \geq 1$, $r_i \geq 2$ and $t_i \geq 1$
\end{corollary}

\begin{remark} {\em It is clear that $\trip$ can be defined for $p=1$ as
well. However, it is not hard to see that $\mathcal{W}(1)$ is the
Virasoro vertex operator algebra $L(1,0)$, while the vertex
algebra generated by $E$,$F$ and $H$ and $\omega$ is all of $V_L$.
In both cases we do not get a new vertex algebra.}
\end{remark}

\section{$C_2$-cofiniteness of $\mathcal{W}(p)$}

As usual, for a vertex operator algebra $V$ we let
$$C_2(V)=\{ a_{-2}b : a, b \in V \}.$$
It is a fairly standard fact (cf. \cite{Zhu}) that $V/C_2(V)$ has
a Poisson algebra structure with the multiplication
$$\bar{a} \cdot \bar{b}=\overline{a_{-1} b},$$
where $-$ denotes the natural projection from $V$ to $V/C_2(V)$. If
${\rm dim}(V/C_2(V))$ is finite-dimensional we say that $V$ is
$C_2$-cofinite.

The aim of this section is the following result
\begin{theorem} \label{c2}
The triplet VOA $\mathcal{W}(p)$ is $C_2$-cofinite.
\end{theorem}
\noindent {\em Proof.} Proposition \ref{stgen} implies that
${\trip}/C_2(\trip)$ is a commutative algebra  generated by the
set $\{ \bar{E}, \bar{F}, \bar{H}, \bar{\omega} \}$. By using
relation (\ref{EFH}) and commutativity of ${\trip} /C_2(\trip)$,
we obtain that
$$\bar{E} ^{2} = \bar{F} ^{2} = 0. $$
Since
$$ Q^{2}(F_{-1} F) = E_{-1} F + F_{-1} E + 2 H_{-1} H =0 $$
we also have that
$$ \bar{H} ^{2} = - \bar{E} \bar{F}$$
which implies that
$$\bar{H} ^{4}=0.$$

Moreover, the description of Zhu's algebra from \cite{A-2003}
implies that
$$ \bar{H} ^{2} = C_p \bar{\omega} ^{2p-1}, \quad (C_p \ne 0).$$
Since $\bar{H} ^{4} = 0$, we conclude that $\bar{\omega}
^{4p-2}=0$. Therefore, every generator of the commutative algebra
${\trip} /C_2(\trip)$ is nilpotent and therefore ${\trip}
/C_2(\trip)$ is finite-dimensional.
 \epfv
\begin{remark}
{\em We should say here that Theorem \ref{c2} was discussed by
Carqueville and Flohr  \cite{CF} by using a different circle of
ideas. Their (rather lengthy) argument was based on analysis of
characters.}
\end{remark}

\section{Zhu's algebra $A(\trip)$ and classification of irreducible $\trip$-modules}

In this part we classify all irreducible $\trip$-modules. Our
approach combines explicit methods (we will eventually realize all
irreducible $\trip$-modules inside the irreducible modules for the
lattice vertex algebra $V_L$) together with an extensive use of
Zhu's associative algebra \cite{Zhu}. Here and throughout the
paper we shall assume some knowledge of the theory of
$C_2$-cofinite vertex algebras, in particular we will use results
from \cite{ABD} and references therein.

Let us recall some fairly standard notation and the definition of
Zhu's algebra $A(V)$ associated to a vertex operator algebra $V$.
As in \cite{AdM} we shall always assume that $$V=\coprod_{ n \in
{\N} } V_n, \ \ \mbox{where} \ \ V_n = \{ a \in V \ \vert \ L(0) a
= n a \}.
$$ For $a \in V_n$, we shall write $\wt (a) = n$ or ${\rm
deg}(a)=n$.

For a homogeneous element $a \in V$ we define the following
bilinear maps $* : V \otimes V \rightarrow V$, $\circ : V \otimes
V \rightarrow V$ as follows: \bea
a*b &:= & \mbox{Res}_x  Y(a,x) \frac{ (1+x) ^{\wt (a) }}{x} b
= \sum_{i = 0} ^{\infty} { \wt (a) \choose i} a_{ i-1} b
, \nonumber \\
a\circ b &: = & \mbox{Res}_x  Y(a,x) \frac{ (1+x) ^{\wt (a) }}{x
^{2}} b
= \sum_{i = 0} ^{\infty} { \wt (a) \choose i} a_{ i-2} b
. \nonumber
 \eea
Next, we extend $*$ and $\circ$ on $V \otimes V$ linearly, and
denote by $O(V)\subset V$ the linear span of elements of the form
$a \circ b$, and by $A(V)$ the quotient space $V / O(V)$. The
space $A(V)$ has a unitary associative algebra structure, the
Zhu's algebra of $V$. The image of $v \in V$, under the natural
map $V \mapsto A(V)$ will be denoted by $[v]$.

Assume that $M= \oplus_{n \in {\N}} M(n)$ is a ${\N}$--graded
$V$--module. Then the top component $M(0)$ of $M$ is a
$A(V)$--module under the action $[a]\mapsto o(a) =a_{\wt (a)-1}$
for homogeneous $a$ in $V$. Moreover, there is ono-to-one
correspondence between irreducible $A(V)$--modules and irreducible
$\N$--graded  $V$--modules (cf. \cite{Zhu}).

Unlike the Poisson algebra $V/C_2(V)$ that inherits the
$V$-grading,  Zhu's associative algebra $A(V)$ is not graded, but
it does admit an increasing filtration. The corresponding
(commutative) associative graded algebra is then denoted by ${\rm
gr}_{\bullet}(A(V))$ (see \cite{Abe}). Then we have a natural
epimorphism of commutative algebras
$$\pi : V/C_2(V) \longrightarrow {\rm gr}_{\bullet} A(V).$$
This and the previous discussion gives the following useful result
(see \cite{Abe} or \cite{Miy} for details).
\begin{proposition} \label{abe}
Let $V$ be strongly generated by the set $S$. Then $A(V)$ is
generated by the set $$\{ [a], a \in S \}.$$ Moreover, if $V$ is
$C_2$-cofinite
$${\rm dim}(V/C_2(V)) \geq {\rm dim}(A(V)).$$
\end{proposition}

Now we specialize $V=\trip$. Our first goal is to gain some
information about the Zhu's algebra (we will have to work more to
obtain additional relations). First, we recall (cf. \cite{Zhu})
that for $a \in V$ homogeneous
%
% promjena z u x
%
\bea \label{rel-z1}\mbox{Res}_x Y(a,x) b \frac{(x+1)^{\wt (a)}}{x
^{2 +n}} \in O({\trip}) \quad \forall a, b \in {\trip}, n \ge 0.
\eea

This implies the following lemma:

\begin{lemma} \label{pomocna-zhu1} We have
$$ Q ^2 e^{-2 \alpha} \in O(\trip) . $$
\end{lemma}
\noindent {\em Proof.} From the relations \bea &&
e^{-\alpha}_{-2p-1}
e^{-\a} = e^{-2 \a}, \nonumber \\
&& e^{-\alpha}_{n} e^{-\a} =0, \ \ \ n \geq -2p, \nonumber \eea and
(\ref{rel-z1}) we obtain
$$e^{-2 \alpha} \in O(\trip).$$
Since $Q$ preserves $O(\trip)$, the proof follows.  \epfv

%
%\bea \label{rel-z2} && Q^{2} e^{-2 \a} = E_{-2 p -1} F + F_{-2p-1} E
%+ 2 H_{-2p-1} H. \eea Also, because of
%$$Q^2(F_{i}F)=0, \ \ i \geq -2p,$$
%we get
%%
%\bea \label{rel-z3} && E_i F + F_i E + 2 H_i H = 0, \quad \forall i
%\ge -2p. \eea By using (\ref{rel-z1})-(\ref{rel-z3}), we have that
%$Q^{2 } e^{-2 \a} \in O(\trip)$. \epfv

By using  \cite{A-2003} and \cite{AdM} one can obtain the
following lemma.

\begin{lemma} \label{pomocna-zhu2}
We have \bea && H_{-2p-1} H \in
\mathcal{U}(Vir).Q^{2} e^{-2 \a} \oplus \mathcal{U}(Vir).{\bf 1} \subset M(1). \nonumber \\
&& H_{-2p-1} H = C Q^{2} e^{-2 \a} + v', \nonumber \eea where $C
\neq 0$ and $v' \in \mathcal{U}(Vir).{\bf 1}$.
\end{lemma}
\noindent {\em Proof.} We already know (cf. \cite{A-2003},
\cite{AdM}), \be \label{rel-hh} H_{-2p-1}H=C Q^{2} e^{-2
\a}+v'+v'', \ee where $C$ is a nonzero complex number, $v'  \in
\mathcal{U}(Vir).{\bf 1}$ and $v'' \in \mathcal{U}(Vir).H$.
(Remember that $Q^2e^{-2 \alpha}=u^{2,2}_0$ is a highest weight
vector for the Virasoro algebra.) Notice that each vector on the
right hand-side of (\ref{rel-hh}) is of conformal weight $6p-2$.
We also recall
$$H_i H \in \mathcal{U}(Vir).{\bf 1}, \ \ i \geq -2p.$$
Suppose that $v'' \neq 0$. The vector $v''$ cannot be singular so
then either $L(1)v''$ or $L(2)v''$ are nontrivial. Suppose that
$L(1).v'' \neq 0$. Then the relation
$$L(1) H_{-2p-1}H=a H_{-2p} H \in \mathcal{U}(Vir).{\bf 1}, \ \ \mbox{for \ some} \ \ a \in
\mathbb{C},$$ together with (\ref{rel-hh}) imply  $L(1)v'' \in
\mathcal{U}(Vir).{\bf 1}$, contradicting $L(1)v'' \in
\mathcal{U}(Vir).H$. The same way we treat the case $L(2)v'' \neq
0$. \epfv

%\begin{lemma}\label{pomocna-zhu2} We have
%%
%\bea \label{rel-z4} && H_{-2p-1} H = C Q^{2} e^{-2 \a} + v', \quad
%\mbox{where} \ C\ne 0, v' \in U(Vir). {\bf 1} \eea
%\end{lemma}
%{\em Proof.} By using  \cite{A-2003} and \cite{AdM} we see that
%%
%$$H_{-2p-1} H = C Q^{2} e^{-2 \a} + v + v'$$
%%
%where $C$ is a non-zero constant, $v \in U(Vir) . H$, $v' \in
%U(Vir) {\bf 1}$. Assume that $v \neq 0$. Then we can take a
%homogeneous  $f \in U(Vir_{\ge 0})$ such that $f v = H$. Since $H$
%is a primary vector we have that $f H_{-2p-1} H = A H(0) H$ for
%certain  constant $A$. Therefore we get that $ H = A H(0) H - f
%v'$, which implies that $Q H = 0$. This is a contradiction and the
%Lemma holds. \epfv

The next lemma will give a very useful binomial identity.
\begin{lemma} \label{identitet-1}
We have the following identity inside $\mathbb{C}[t]$: \be
\label{t-ind} \widetilde{{\Phi}_p}(t) = \sum_{i=0} ^{2p} (-1) ^{i}
{2p \choose i} {t \choose 4p-1-i} { t \choose 2p + i -1} = A_p {t
\choose 3p-1} {t+p \choose 3p-1}, \ee where
$A_p = (-1) ^{p} \frac{ {2p  \choose p} }{{4p-1 \choose p}}$.
\end{lemma}
\no {\em Proof.} Let $g_p(t) = A_p {t \choose 3p-1} {t+p \choose
3p-1}$. It is clear that both $g_p(t)$ and
$\widetilde{{\Phi}_p}(t)$ are inside $\mathbb{C}[t]$. First we
notice that for $t= 3p-1$ there is only one nontrivial term in
(\ref{t-ind}), so
\bea \label{3p-1} && \widetilde{{\Phi}_p}(3p-1) = (-1) ^{p} { 2p
\choose p} = A_p {4p-1 \choose 3p-1}= g_p(3p-1).\eea
By using straightforward calculation we get the following
recursion:

\bea \label{rek-1} &&
(t+1) (p + t+1) \widetilde{{\Phi}_p}(t) = (2p-t-2) (3p-t-2)
\widetilde{{\Phi}_p}(t+1).
 \eea
By using (\ref{3p-1}) and the fact that  $g_p(t)$ also satisfies
the same recursion we conclude that $\widetilde{{\Phi}_p} (t) =
g_p(t)$ for infinitely values of $t \in {\C}$. The proof follows.
\epfv
\begin{theorem} \label{zhu1}
In Zhu's algebra $A(\trip)$ we have the following
relation
$$f_p([\omega])=0,$$
where \be \label{relzhu1} f_p(x)=\prod_{i=0}^{3p-2}(x-h_{i+1,1}).
\ee
\end{theorem}
\no {\em Proof.}
 By using   Lemma \ref{pomocna-zhu2} we have that (in $A(\mip)$, where $\mip$ is the singlet algebra)
 \bea
 \label{rel-z5} && [Q^{2} e^{-2 \a}] = {\Phi}_p([\omega])
 \eea
for certain ${\Phi}_p \in {\C}[x]$, $\deg {\Phi}_p \le 3 p-1$.
We will see that ${\Phi}_p$ is an non-trivial polynomial of degree
$ 3p-1$, and find explicit formulae for this polynomial. Since
$Q^{2}e^{-2\a} \in \mip \subset M(1)$,
 we shall evaluate the action of $Q^{2} e^{-2 \a}$ on top levels of $\overline{M(1)_p}$--modules $M(1,\l)$. Let $v_{\l}$ be the highest weight vector in
 $M(1,\l)$. We use the method of \cite{A1} and  \cite{A-2003}. By using (\ref{prod-vo}) and direct calculation we see that
 \bea && o( Q^{2} e^{-2\a}) v_{\l} =  {\Phi}_p (\tfrac{1}{4p} (t ^{2} - 2 t (p-1) ) v_{\l}\nonumber \\
&& =\mbox{Res}_{z_1} \mbox{Res}_{z_2}  \ (z_1 z_2) ^{-4p} (z_1 -
z_2) ^{2p} (1+z_1) ^{t} (1+z_2) ^{t} v_{\l} = \widetilde{{\Phi}_p}(t) v_{\l} \nonumber \\
&& \mbox{where} \quad \ \widetilde{{\Phi}_p}(t) = \sum_{i=0} ^{2p}
(-1) ^{i} {2p \choose i} {t \choose 4p-1-i} { t \choose 2p + i -1},
\quad t= \la \l, \a \ra . \nonumber \eea
By using Lemma \ref{identitet-1} we get that
 $$  {\Phi}_p (\tfrac{1}{4p} (t ^{2} - 2 t (p-1) ) =\widetilde{{\Phi}_p}(t) = B_p f_p ( \tfrac{1}{4p} (t ^{2} - 2 t (p-1) ), $$
 where the nonzero constant $$B_p=  (-1) ^{p} \frac{ {2p \choose p} (4p) ^{3p-1} }{{4p -1 \choose p} (3p-1) !
 ^{2}}.$$
 This proves that ${\Phi}_p$ is a non-trivial polynomial of degree
$3p-1$ and that  in $A(\mip)$ we have
 \bea \label{rel-q2}  && [Q^{2} e^{-2\a}] ={\Phi}_p ([\omega]) =  B_p
 f_p([\omega]).\eea
Since $O(\mip) \subset O(\trip)$, the embedding $\mip \subset
\trip$ induces an algebra homomorphism $A(\mip) \rightarrow
A(\trip)$. Applying this homomorphism to relation (\ref{rel-q2})
and using Lemma \ref{pomocna-zhu1} we get that $f_p([\omega]) = 0$
in $A(\trip)$.
 \epfv
%
% The statement of theorem now follows from Lemma
%\ref{pomocna-zhu1}. \epfv

\vskip 5mm

\begin{corollary} \label{bound} The relation
$$\bar{\omega}^{3p-1}=0$$
holds inside $\trip/C_2(\trip)$.
\end{corollary}
\noindent
\no {\em Proof.} By using Lemma \ref{pomocna-zhu2} and the proof
of Theorem \ref{zhu1} we obtain
\bea  \label{rel-c2} && Q^{2} e^{-2\a} = A H_{-2p-1} H + \left(B_p
L(-2) ^{3p-1} + w \right).{\bf 1} \eea
where $A\ne 0$, $B_p \ne  0$ and $w$ is a linear combination of
monomials in $\mathcal{U}(Vir_{\leq -2})$ each containing at least
one generator $L(-n)$, $n \geq 3$. Clearly $w. {\bf 1} \in
C_2(\trip)$. Since
$$H_{-2p-1} H \in C_{2}(\trip)$$
and
$$ Q^{2} e^{-2\a} = E_{-2p-1} F + F_{-2p-1} E + 2 H_{-2p-1}
H \in C_{2} (\trip),$$
the relation (\ref{rel-c2}) implies  $\bar{\omega} ^{3p-1} = 0 $.
\epfv

%{\em Proof.} Since ${\rm deg}(f_p)=3p-1$ and ${\rm
%deg}(e^{-2 \alpha})=6p-2$, Theorem \ref{zhu1} implies that
%$$\left(L(-2)^{3p-1}+ w \right).{\bf 1} \in O(\trip),$$
%where $w$ is a linear combination of monomials in
%$\mathcal{U}(Vir_{\leq -2})$ each containing at least one
%generator $L(-n)$, $n \geq 3$. Now, the discussion prior to
%Proposition \ref{abe} implies $\bar{\omega}^{3p-1}=0$ \epfv

\begin{remark} {\em Notice that Corollary \ref{bound} gives a better upper bound on
${\rm dim}(\trip/C_2(\trip))$ compared to the one that could be
extracted from Theorem \ref{c2}.}
\end{remark}

The polynomial (\ref{relzhu1}) can be rewritten more conveniently
as \be \label{zhufp}
f_p(x)=(x-h_{p,1})(x-h_{2p,1})\prod_{i=1}^{p-1} (x-h_{i,1})^2
\prod_{i=2p+1}^{3p-1} (x-h_{i,1}). \ee Already from (\ref{zhufp})
and $f_p([\omega])=0$ it seems feasible to expect some non
$L(0)$-diagonalizable weak $\trip$-modules entering the picture.
More precisely, notice that $f_p([\omega])$ will annihilate not
only every top component of a $\trip$-module with lowest conformal
weight $h_{i,1}$, $1 \leq i \leq 3p-1$, but also any lowest weight
subspace of {\em generalized} conformal weight $h_{i,1}$, $1 \leq
i \leq p-1$ with Jordan blocks--with respect to $L(0)$--of size at
most two (cf. second powers in (\ref{zhufp})). We will get back to
this problem in the later sections.

The previous theorem gives also an important information about the
possible lowest conformal weights of irreducible $\trip$-modules.
Now, we construct a family of $\trip$-modules with lowest
conformal weights being precisely the roots of $f_p(x)$.

Recall (\ref{ga-def}). Since $\trip \subset V_L$, every
(irreducible) $V_L$-module is of course a $\trip$-module. Thus it
is natural to examine the structure of $V_{L+ \ga_i}$, viewed as a
$\trip$-module. By using the formula
$${\rm deg}(e^{m \alpha})=m^2p-m(p-1),$$
which holds for every $m \in \frac{\mathbb{Z}}{2p}$, we see that
$${\rm deg}(e^{\ga_i})=\frac{i(i-2p+2)}{4p}= h_{i+1,1}$$
and
$${\rm deg}(e^{\ga_{2p-2-i} })=\frac{i(i-2p+2)}{4p}= h_{i+1,1}.$$
From the last two formulas it is not hard to see that
$$V_{L+\ga_i}, \ \ i \neq 2p-1,$$
have a one-dimensional lowest weight subspace spanned by
$e^{\ga_i}$, while
$$V_{L+\ga_{2p-1}}$$
has a two-dimensional top subspace spanned by $e^{(2p-1)
\alpha/2p}$ and $e^{-\alpha/2p}$. We also have an  isomorphism \be
\label{dual-iso} V_{L+\ga_i}' \cong V_{L+\ga_{2p-2-i} }, \ee where
$'$ stands for the contragradient dual. Thus the only irreducible
self-dual $V_L$-module is $V_{L+\ga_{p-1}}$.

Now, we identify certain submodules and subquotients of each
$V_{L+ \ga_{i-1}}$.

For $i=1,..,p-1$ we set
\bea && \Lambda(i):= \overline{V_{L+\ga_{i-1}}} \subset
V_{L+\ga_{i-1}} \nonumber \\
&&\Pi(p-i) : = \overline{V_{L + \ga_{2p-i-1}} }\subset V_{L +
\ga_{2p-i-1} } ,\nonumber \eea
and
\bea
&&\Lambda(p) =  V_{L+\ga_{p-1}}, \nonumber \\
&& \Pi(p) = V_{L+\ga_{2p-1}}.\nonumber
\eea
% the $\trip$-module generated by
%$e^{(i-1) \alpha/2p}$
(Here we are actually using the notation introduced in
\cite{FHST}). Thus we obtain the following short exact sequences
of vector spaces:
$$0 \longrightarrow \Lambda(i) \longrightarrow V_{L+\ga_{i-1}}
\longrightarrow  V_{L+\ga_{i-1}}/\Lambda(i) \longrightarrow 0 ,$$
$$0 \longrightarrow \Pi(p-i) \longrightarrow V_{L+\ga_{2p-i-1}}
\longrightarrow  V_{L+\ga_{2p-i-1} }/\Pi(p-i) \longrightarrow 0.$$
It is important to notice that the lowest conformal weight of
$\Pi(p-i)$ is
$$h_{2p+i,1} =h_{2p-i,3}=\frac{(i+1)(2p+i-1)}{4p}.$$

\vskip 5mm

\begin{theorem} \label{ireducibilni-moduli} For every $1 \leq i \leq p$, both $\Lambda(i)$ and $\Pi(i)$ are irreducible
$\trip$-modules.
%For $1 \leq i \leq p-1$ we also have short exact sequences \bea
%\label{short1} && 0 \longrightarrow \Lambda(i) \longrightarrow
%V_{L+(i-1)\alpha/2p} \longrightarrow  \Pi(p-i) \longrightarrow 0
%\\ && \label{short2} 0 \longrightarrow \Pi(p-i) \longrightarrow V_{L+(2p-i-1)
%\alpha/2p} \longrightarrow  \Lambda(i) \longrightarrow 0. \eea
\end{theorem}
\noindent {\em Proof.} Assume  that $0\le i \le p-1$.
%
% pocetak dokaza je isti kao u AM-2007
%
In \cite{AdM} we argue that the space of intertwining operators
$$I \ { L(c_{p,1},h) \choose  L(c_{p,1},2 p -1) \ \  L(c_{p,1},
h_{i+1,2n+1})}$$ is non-trivial if and only if $h=h_{i+1,2n-1}$,
$h=h_{i+1,2n+1}$ or $h=h_{i+1,2n+3}$, for $n \geq 1$. Since the
multiplicities of these fusion rules are always one, we may write
formally:
\bea  \label{fusion-rules} &&
L(c_{p,1},2 p -1) \times L(c_{p,1}, h_{i+1,2n+1} ) =\nonumber \\
 && L(c_{p,1}, h_{i+1,2n-1} ) \oplus L(c_{p,1}, h_{i+1,2n+ 1} ) \oplus L(c_{p,1}, h_{i+1,2n+3} ) \quad (n \ge 1).\label{f-r}\eea
The same results has been known by physicists (cf. \cite{F1}). The
fusion rules (\ref{f-r}) implies that
$$X_j \Lambda(i+1) \subset \Lambda(i+1) , \quad
\mbox{for every}  \ j \in {\Z}, \quad \mbox{where} \ X=E, F \
\mbox{or} \ H.$$ Since $\trip$ is generated by $\omega$, $E$, $F$
and  $H$ we have that $\Lambda(i+1)$ is an $\trip$--module. In
order to prove that $\Lambda(i+1)$ is irreducible, we shall first
prove that $e^{\ga_i}$ is a cyclic vector in $\Lambda(i+1)$, i.e.,
$\Lambda(i+1) = \langle e^{\ga_i} \rangle$. Assume that $\langle
e^{\ga_i}  \rangle \ne \Lambda(i+1) $. This implies that $M =
\Lambda(i+1) / \langle e^{\ga_i} \rangle$ is a non-trivial
${\N}$--graded ${\trip}$--module. By Theorem \ref{zhu1} it follows
that $M = \oplus_{ n \in {\N} } M(h'+n)$, where $h' = h_{j+1,1}$
for $0 \le j \le 3p-2$. Recall that, as a module for the Virasoro
algebra, $\Lambda(i+1)$ is generated by the family of singular
vectors of weights $h_{i+1,2n+1}$, $n \in {\N}$.  Therefore every
non-trivial vector from the top level $\bar{v} \in M(h')$ has the
form $\bar{v} = v + \langle e^{\ga_i} \rangle$, where $v$ is a
singular vector for the Virasoro algebra of weight $h'$. This
leads to a contradiction since
$ h' = h_{i+1,2n+1}$ if and only if  $n = 0$, and every vector of
weight $h_{i+1,1}$ in $\Lambda(i+1) $ is proportional to $e^{\ga_i}
\in \langle e^{\ga_i} \rangle$.
 Thus $\Lambda(i+1)= \langle e^{\ga_i} \rangle$.

Assume that $N$ is any non-trivial ${\trip}$--submodule of
$\Lambda(i+1)$. Then $N$ is also a $ {\N}$--graded
${\trip}$--module. Again, by Theorem \ref{zhu1}, $N= \oplus_{ n \in
{\N}} N(h'+n)$, where $h' = h_{j+1,1}$ for $0 \le j \le 3p-2$. The
same arguments as before show that $h'=h_{i+1,1}$, which implies
that $e^{\ga_i} \in N$ and therefore $N=\Lambda(i)$.

By using  the fusion rules, one can also prove that $\Pi(p-i)$ is
a $\trip$--module, where $0 \le i \le p-1$. The top component
$\Pi(p-i)(0)$ is an irreducible two dimensional $A(\trip)$--module
with conformal weight $h_{2p + i,1}=h_{2 p-i,3}$. Theorem
\ref{str-ff-2} implies that as a module for the Virasoro algebra
$\Pi(p-i)$ is generated by the family of singular vectors of
weights $h_{2 p-i,2n+1}$, $n \ge 1$.

Let $j \in \{1, \dots, 3p-1\}$. By using the fact that
$$ h_{2 p-i, 2n +1} = h_{j,1} \quad \mbox{iff} \quad n=1, j =
2p +i ,$$
and a completely analogous proof as in the case of $\Lambda(i+1)$
one can prove that $\Pi(p-i)$ is an irreducible $\trip$--module.
%
%Similarly we see that $\Pi(i+1)$ is an irreducible $\trip$--module.
%Then (\ref{short1}) and (\ref{short2}) follow from Theorems
%\ref{str-ff-1} and \ref{str-ff-2}.
\epfv

Applying the previous theorem in the case of $\trip =\Lambda(1)$
we get:
\begin{corollary}
The vertex operator algebra $\trip$ is simple.
\end{corollary}
%The module $\Lambda(i)$ is the maximal semisimple submodule of
%$V_{L+i \alpha/2p}$, viewed as Virasoro algebra module (We should
%give a proof here). From the Virasoro algebra structure we can
%easily see that $\Pi(s)$

The next result will be proven in Theorem \ref{zhu-relations}.
\begin{proposition} \label{zhu-komut}
In Zhu's associative algebra we have \bea &&
[H]*[F]-[F]*[H]=-2q([\omega])[F], \\
&& [H]*[E]-[E]*[H]=2q([\omega])[E] \\
&& [E]*[F]-[F]*[E]=-2q([\omega])[H]. \eea where $q$ is a polynomial
of degree at most $p-1$. \end{proposition}

It is known that the top component $M(0)$ of a $\mathbb{Z}_{\geq
0}$-gradable $\trip$-module $M$ carries a structure of
$A(\trip)$-module. In particular, the top component $\Pi(i)(0)$
and $\Lambda(i)(0)$ are $A(\trip)$-modules. More precisely, the
top component of $\Lambda(i)$ is $1$-dimensional and has conformal
weight $h_{i,1}$ for $i=1, \dots, p$.
 On the other
hand, the top component  of $\Pi(i)$ is $2$--dimensional and has
conformal weight $h_{p+i,3}=h_{3p-1-i,1}$ for $i=1, \dots, p$. So
we  get:

\begin{proposition} \label{konstr-ired}
\item[(1)]For every $ 1 \le i \le p $ the top component
$\Lambda(i)(0)$ of $\Lambda(i)$ has lowest conformal weight
$h_{i,1}$. Moreover $\Lambda(i) (0)$  is an
  irreducible $1$--dimensional $A({\trip})$--module spanned by
  the  vector $e^{\ga_{i-1}}$.

\item[(2)] For every $2p \le i\le 3p-1$ the top component
$\Pi(3p-i) (0)$ of $\Pi(3p-i)$ has lowest conformal weight
$h_{i,1}$. Moreover, $\Pi(3p-i)(0)$  is an irreducible
$2$--dimensional $A({\trip})$--module spanned by   the vectors
$ e^{\ga_{2p-1-i}}$ and $ Q e^{\ga_{2p-1-i}}$.
\end{proposition}

\begin{theorem} \label{class-ired} The set $$\{ \Pi(i)(0) : 1 \leq i \leq p  \} \cup
\{ \Lambda(p-i)(0): 1 \leq i \leq p \}$$ provides, up to
isomorphism, all irreducible modules for  Zhu's algebra
$A(\trip)$.
\end{theorem}

\no {\em Proof.} Assume that $U$ is an irreducible
$A(\trip)$--module. Relation $f_p ([\omega]) = 0$ in $A(\trip)$
implies that
$$L(0) \vert U = h_{i,1} \ \mbox{Id}, \quad \mbox{for} \quad i \in \{1, \dots, p \} \cup \{2p, \dots, 3p-1\}.$$

Assume first that $ 2p \le i \le 3p-1$. By combining Propositions
\ref{zhu-komut} and \ref{konstr-ired} we have that $q(h_{i,1}) \ne
0$. Define
$$ e= \frac{ 1}{ \sqrt{2} q(h_{i,1})} E, \quad f= -\frac{1}{ \sqrt{2}q(h_{i,1})} F, \quad h= \frac{1}{q(h_{i,1})} H .$$

Therefore $U$ carries the structure of an irreducible,
$\goth{sl}_2$--module with the property that $e^2 = f^2 = 0$  and
$h \ne 0$ on $U$. This easily implies that $U$ is $2$--dimensional
irreducible $\goth{sl}_2$--module. Moreover,  as an
$A({\trip})$--module $U$ is isomorphic to  $\Pi(3p-i)(0)$.

Assume next that $ 1 \le i \le p$. If $q(h_{i,1}) \ne 0$, as above
 we conclude that $U$ is an irreducible $1$--dimensional $sl_2$--module. Therefore $U \cong \Lambda(i)(0)$.

 If $q(h_{i,1})=0$
 from Proposition \ref{zhu-komut} we have that the action of generators of $A({\trip})$ commute on $U$. Irreducibility of
  $U$ implies that $U$ is $1$-dimensional. Since $[H], [E] ^2, [F] ^2$ must act trivially on $U$, we conclude that $[H], [E] , [F] $ also act trivially on $U$. Therefore
  $U \cong \Lambda(i)(0)$.
  \qed

\begin{remark}
{\rm In Theorem \ref{zhu-relations} below we shall see that in
fact $q(h_{i,1}) \ne
 0$ for every $i \in \{1, \dots, p\}$, which can be used to give a shorter proof of
Theorem \ref{class-ired}.}
\end{remark}

Applying Zhu's theory \cite{Zhu} on Theorem \ref{class-ired} and
using irreducibility result from Theorem \ref{ireducibilni-moduli}
we get the classification of all irreducible $\trip$--modules (the
same result was stated in \cite{FHST}).

\begin{theorem} \label{class-ired-modules} The set $$\{ \Pi(i) : 1 \leq i \leq p  \} \cup
\{ \Lambda(p-i)  : 1 \leq i \leq p \}$$ provides, up to
isomorphism, all irreducible modules for the vertex operator
algebra $\trip$.
\end{theorem}

\section{A description of the category of ordinary $\trip$-modules}

In the previous section we classified simple objects in the
category of $\trip$-modules. Here we derive additional results
about reducible modules. In \cite{AdM} we classified irreducible
self-dual $\overline{M(1)_p}$-modules. For the triplet vertex
algebra we have a very different result.

\begin{corollary} All irreducible $\trip$-modules are self-dual.
\end{corollary}
\noindent {\em Proof.} It is sufficient to show that $\Pi(i)'(0)
\cong \Pi(i)(0)$ and $\Lambda(i)'(0) \cong \Lambda(i)(0)$, where
$W'$ denotes the dual module of $W$. Both isomorphisms follow
directly Theorem \ref{ireducibilni-moduli} and $h_{i,1} \neq
h_{j,1}$, $i \neq j$, $i,j \in \{ 1 \leq k \leq p  \} \cup \{ 2p
\leq k \leq 3p-1 \}$. \epfv

\begin{proposition} \label{ses}
For $1 \leq i \leq p-1$ we have the following non-split short
exact sequences of $\trip$-modules: \bea \label{short1} && 0
\longrightarrow \Lambda(i) \longrightarrow V_{L+(i-1)\alpha/2p}
\longrightarrow \Pi(p-i) \longrightarrow 0
\\ && \label{short2} 0 \longrightarrow \Pi(p-i) \longrightarrow V_{L+(2p-i-1)
\alpha/2p} \longrightarrow  \Lambda(i) \longrightarrow 0. \eea
\end{proposition}
\noindent {\em Proof.} From Theorem \ref{ireducibilni-moduli} we
know that both $\Lambda(i)$ and $\Pi(i)$ are irreducible
$\trip$-modules. Consider $\Psi(p-i):=V_{L+\ga_{i-1}}/\Lambda(i)$.
It is clear (cf. Theorems \ref{str-ff-1} and \ref{str-ff-2}) that
as a Virasoro algebra module $\Psi(p-i)$ is isomorphic to
$\Pi(p-i)$. Also, as $A(\trip)$-modules $\Psi(p-i)(0) \cong
\Pi(p-i)(0)$. Now, irreducibility of $\Pi(p-i)$ gives $\Psi(p-i)
\cong \Pi(p-i)$. Thus, we have (\ref{short1}). Similarly, we show
(\ref{short2}). From Theorem \ref{str-ff-1} and Theorem
\ref{str-ff-2} it is clear that the sequences are non-split. \epfv

The previous proposition seems to be know as well as the next
result (cf. \cite{FHST}).

%\begin{proposition} The category of ordinary (non-logarithmic)
%$\trip$-modules consists of $p+1$ linkage classes $\mathcal{F}_i$,
%$i=0,...,p$, where $\Lambda(p)$ and $\Pi(p)$ are the only
%irreducible objects in $\mathcal{F}_0$ and $\mathcal{F}_{p}$,
%respectively, while for $1 \leq i \leq p-1$ each linkage class
%$\mathcal{F}_i$ contains $\Lambda(i)$ and $\Pi(p-i)$.
%\end{proposition}

Now, we recall several standard facts about category $\mathcal{O}_c$
for the Virasoro algebra. A Virasoro module $M$ is said to be an
object in $\mathcal{O}_c$ if
\begin{itemize} \item[(i)] $M=\oplus_{n \in \mathbb{C}} M_n$, \
$M_n=\{ v \in M: L(0)v=nv \}$, \ ${\rm dim}(M_n)< + \infty$,
\item[(ii)] The central element acts as multiplication by $c$,
\item[(iii)] There exist $\lambda_1,...,\lambda_k \in \mathbb{C}$
such that if $M_n \neq 0$ then $n \geq \lambda_i$ for some $i$.
\end{itemize}

We denote by $M(c,h)$ the Verma Virasoro module with lowest
conformal weight $h$ and central charge $c$. We introduce a
partial ordering $"\preceq"$ on the set of weights
($=\mathbb{C}$). We say that $h' \preceq h$ if
$$L(c,h') \ \ \mbox{is a
subquotient of} \ \  M(c,h).$$ Extend the partial ordering
$"\preceq"$ to an equivalence relation $\sim$ on the set of
weights.
% Then $"\preceq"$ is an equivalence relation on the set of
%weights.
Let $[h]$ denote an equivalence class and by
$\mathcal{F}$ the set of all equivalence classes. Then it is known
that every $W \in Obj(\mathcal{O}_c)$ admits {\em block
decomposition}:
$$W= \bigoplus_{[h] \in \mathcal{F}} W_{[h]}, \ \ W_{[h]} \in \mathcal{O}^{[h]}_c,$$
where $M \in \mathcal{O}^{[h]}_c$ if every irreducible subquotient
$L(c,h')$ of $M$ satisfies $h' \in [h]$. Each
$\mathcal{O}^{[h]}_c$ is in fact a full subcategory of
$\mathcal{O}_c$.

We would like to obtain a similar description for the category of
{\em ordinary} (i.e., non-logarithmic) $\trip$-modules. Since
$\trip$ is not a Lie algebra with triangular decomposition one
cannot define $\trip$-blocks as above by using Verma modules.

\begin{lemma} Every ordinary $\trip$-module, viewed as a Virasoro algebra module, is
an object in $\mathcal{O}_{c_{p,1}}$.
\end{lemma}
\noindent {\em Proof.} We only have to check that there exists
$\lambda_1$,...,$\lambda_k$ but this is clear since $\trip$ is
$C_2$-cofinite. \epfv

The previous lemma indicates that we should first try to describe
Virasoro blocks for irreducible $\trip$-modules.

The first important observation we make is that each irreducible
$\trip$-module belongs to a unique (Virasoro) block. This is a
consequence of Theorem \ref{str-ff-1} and  \ref{str-ff-2}, and
the fact that
%
% dodano
%
\bea && \label{nprec } h_{i,1} \notin [h_{j,1}], \quad \mbox{for}
\ i, j \in \{1,\dots, p\} \cup \{2p\}, \ i \ne j. \eea

 Actually we can say more; $\Lambda(p)$ and $\Pi(p)$
live in two distinct blocks and the remaining irreducible modules
are distributed in additional $p-1$ blocks, such that $\Lambda(i)$
and $\Pi(p-i)$ are in the same block. These $p+1$ blocks are
represented by $h_{i,1}$, $i=1,...,p-1$, $h_{p,1}$ and $h_{2p,1}$.
Our next result says that the same decomposition persists at the
level of $\trip$-modules.

\begin{theorem} The category of ordinary $\trip$-modules contains
precisely $p+1$ blocks, so that every $\trip$-module $W$
decomposes as a direct sum of $\trip$-modules: \be \label{block}
W= W_{[h_{1,1}]} \oplus \cdots \oplus W_{[h_{p,1}]} \oplus
W_{[h_{2p,1}]}, \ee where
$$W_{[h_{i,1}]} \in Obj(\mathcal{O}^{[h_{i,1}]}_{c_{p,1}}).$$
\end{theorem}
\noindent {\em Proof.} We certainly have a decomposition of $W$ on
the level of Virasoro modules (with possibly additional summands).
In view of the previous discussion it suffices to show that $F$,
$H$ and $E$ preserve the summands ${W}_{[h_{i,1}]}$ and that there
are no additional blocks appearing in (\ref{block}). Since $E$,
$F$ and $H$ are primary fields of lowest conformal weight $2p-1$
the action of $\trip$ on $W$ defines a $L(c_{p,1},0)$-intertwining
operator of type
$${ W \choose L(c_{p,1},2p-1) \ \ W_{[h_{i,1}]} }.$$
If $W$ admits a nontrivial summand from a different block it would
follow that there is a nontrivial intertwining operator of the
type
$${ W_{[h_{j,1}]} \choose L(c_{p,1},2p-1) \ \ W_{[h_{i,1}]} }, \ \ h_{j,1} \neq h_{i,1}.$$
But this will contradict to the fusion rules formula
(\ref{fusion-rules}). It follows that $E$, $F$ and $H$ preserve
each block. To show that there are no additional blocks one  uses
the classification of irreducible $\trip$--modules from Theorem
\ref{class-ired-modules}. \epfv
%
%To show that there are no additional blocks one again uses formula
%(\ref{fusion-rules}). \epfv

Since each Virasoro block contains at least one irreducible
$\trip$-modules and there is a nontrivial extension of
$\Lambda(i)$ and $\Pi(p-i)$ inside the single block, (\ref{block})
is the proper block decomposition of a $\trip$-module.
Correspondingly, the category of $\trip$-modules is a block
preserving, subcategory of $\mathcal{O}_{c_{p,1}}$.

\begin{corollary}
For every $p \geq 2$ and $i \neq j$ we have
$${\rm Ext}^1_{\trip}(\Lambda(i),\Lambda(j))={\rm
Ext}^1_{\trip}(\Pi(i),\Pi(j))=0,$$
$${\rm Ext}^1_{\trip}(\Lambda(i),\Pi(p-j))={\rm
Ext}^1_{\trip}(\Pi(i),\Lambda(p-j))=0,$$ where the ${\rm
Ext}$-groups are computed inside the category of ordinary
$\trip$-modules.
\end{corollary}
This result has been stated in \cite{FHST}.

\section{The structure of  $A(\trip)$}

In this section we derive additional relations in Zhu's algebra
needed for a better description of $A(\trip)$. Some relations
obtained here could be also used for classification of irreducible
modules in a simpler fashion.

The following important theorem gives a fairly explicit
description of $A(\trip)$ in terms of generators and relations.
\begin{theorem} \label{zhu-relations}
The Zhu's associative algebra is generated by $[\omega]$, $[F]$,
$[H]$, and $[E]$. Also, we have the following relations in
$A(\trip)$:
\begin{itemize}
\item[(i)]   $[E]^2=[F]^2=0$

\item[(ii)] $[H]^2=C_p P([\omega])$, where
$$P(x)=(x-h_{p,1})\prod_{i=1}^{p-1}(x-h_{i,1})^2 \in \mathbb{C}[x]$$
and $C_p$ is a nonzero rational number.

\item[(iii)] $$[H]*[F]=-[F]*[H]=-q([\omega])*[F],$$
$$[H]*[E]=-[E]*[H]=q([\omega])*[E],$$ where $q(x)$ is a nonzero
polynomial of degree $\leq p-1$ and $$q(h_{i,1}) \neq 0, \ \ 1
\leq i \leq p.$$

\item[(iv)] \bea &&
[H]*[F]-[F]*[H]=-2q([\omega])[F], \nonumber \\
&& [H]*[E]-[E]*[H]=2q([\omega])[E], \nonumber \\
&& [E]*[F]-[F]*[E]=-2q([\omega])[H], \nonumber \eea where $q(x)$ is
as in (iii).

\item[(v)] $$ \prod_{i=2p}^{3p-1}([\omega]-h_{i,1})*[X]=0, \ \ X
\in \{ E,F,H \}.$$

\end{itemize}
\end{theorem}

\noindent {\em Proof.} By Proposition \ref{abe} it is clear that
$A(\trip)$ is generated by $[\omega]$, $[F]$, $[H]$ and $[E]$.

Part (i) is trivial to show. Part (ii) has been established
earlier in \cite{A-2003}. It remains to show (iii), (iv) and (v).

Observe that the screening $Q$ preserves $O(\trip)$, hence it acts
as a derivation of $A(\trip)$. We will abuse the notation and use
the same letter for the projection of $Q$ on $A(\trip)$, so in
particular
$$Q[F]=[H], \ \ Q[H]=[E].$$
Then
$$0=Q([X]*[X])=[H]*[X]+[X]*[H], \ \ X \in \{E,F \},$$
which yields two (easy) assertions in (iii). If we can show that
\be \label{hf} [H]*[F]=-q([\omega])*[F], \ee
 then the formulas
$$Q^2([H]*[F])=[H]*[E]+2 [E]*[H]=-[H]*[E]$$
$$Q^2(q([\omega]*[F])=q([\omega])*[E]$$
yield
$$[H]*[E]=q([\omega])*[E].$$
The proof of relations (\ref{hf}) with $q(x)$ satisfying the
properties as in the theorem is rather technical and it is given
in Appendix at the end of the paper.

For (iv), the first two relations follow directly from (iii), and
for the third observe that
$$Q([H]*[F]-[F]*[H])=[E]*[F]-[F]*[E]=-2q([\omega])*[H].$$
The proof follows.

Next we prove part (v). By part (iii), $q(h_{i,1}) \neq 0$, for
$i=1,..,p$. Hence $q([\omega])$ is a unit in $A(\trip)$. Since the
$Q$-screening preserve $O(\trip) \subset \trip$, and commutes with
the action of the Virasoro algebra, the relation
$$l([\omega])*[F]=0, \ \ l([\omega])=\prod_{i=2p}^{3p-1}([\omega]-h_{i,1})$$
implies
$$l([\omega])*[X]=0, \ \ X \in \{E,F,H \}.$$
Suppose on the contrary that
$$l([\omega])*[F] \neq 0.$$
Since $$0=Q^2 ([F]^2) = [E] * [F] + [F]* [E] + 2 [H]^2,$$ by using
(ii) we get \be \label{efanti} [E]*[F]+[F]*[E]=-2 C_p P([\omega]).
\ee The last assertion gives \be \label{fe}
[F]*[E]=q([\omega])*[H]-C_p P([\omega]). \ee Thus, we have
 \be \label{zhu-v-rel}
l([\omega])*[F]*[E]=l([\omega])(q([\omega])*[H]-C_p
P([\omega]))=l([\omega])*q([\omega])*[H], \ee where in the last
equality we used $$l([\omega])*C_p P([\omega])=C_p
f_p([\omega])=0,$$ which holds in $A(\trip)$. After we multiply
(\ref{zhu-v-rel}) with $[F]$, from the left and use $[F]^2=0$, and
Theorem \ref{zhu-relations} (iii), we obtain
$$0=l([\omega])*q([\omega])*[F]*[H]=l([\omega])*q([\omega])^2*[F].$$
But since $q([\omega])^2$ is a unit in $A(\trip)$ we obtain
$$l([\omega])*[F]=0,$$
a contradiction.

\epfv

\begin{corollary}
Zhu's algebra $A(\trip)$ contains a Lie subalgebra isomorphic to
$\goth{sl}_2$. Correspondingly, any $A(\trip)$--module is
naturally an $\goth{sl}_2$--module.
\end{corollary}
\no {\em Proof.} Since $q(h_{i+1,1}) \ne 0$ for every $i$,  $0 \le
i \le 3p-2$, we have that $q(x)$ is relatively prime with $f_p(x)$
and therefore $q([\omega])$ is an unit in $A(\trip)$. Define
nonzero vectors
$$ e= \frac{1}{ \sqrt{2}} q([\omega]) ^{-1} E, \quad f= -\frac{1}{ \sqrt{2} } q([\omega]) ^{-1} F, \quad h= q([\omega]) ^{-1} H .$$
It is easy to see that $[e,f]=h, [h,f]=-2f$ and $[h,e]=2e$ holds.
Thus ${\rm span}\{e,f,h \}$ is isomorphic to $\goth{sl}_2$. This
also implies that any $A(\trip)$--module, in particular $A(\trip)$
itself, is a $\goth{sl}_2$--module. \epfv

\begin{corollary}
The associative algebra $A(\trip)$ is spanned by
$$\{ [\omega]^i, \ 0 \leq  i  \leq  3p-2 \} \cup \{[\omega]^i*[X], 0 \leq i
\leq p-1, \ X=E,F \ \ {\rm or} \ H \}.$$ Thus, $A(\trip)$ is at most
$(6p-1)$-dimensional
\end{corollary}

Our goal is to describe $A(\trip)$ as a sum of ideals. For these
purposes, for $i=1,...,p-1$, let
$$v_i=(\lambda_i[\omega]+\nu_i)\prod_{j=1,...,3p-1; j \neq i, j \neq
2p-i}([\omega]-h_{j,1}),$$ where $\lambda_i$ and $\nu_i$ are
unspecified constants and
$$w_i=\prod_{j=1,..,3p-1; j \neq i} ([\omega]-h_{j,1}).$$
$$v_i * w_i=w_i * v_i=w_i.$$
We clearly have
$$v_i * v_j=0, \ i \neq j, $$
{ and}
$$w_i^2=0, \ i=1,..,p-1.$$

\begin{lemma} For every $i=1,...,p-1$ there are constants $\lambda_i$, $\nu_i$ such that
$$v_i=(\lambda_i[\omega]+\nu_i)\prod_{j=1,...,3p-1; j \neq i, j \neq
2p-i}([\omega]-h_{j,1}) \neq 0$$ satisfy
$$v_i * v_j = \delta_{i,j} v_i,$$
and
$$w_i * v_j=v_j * w_i = \delta_{i,j} w_i,$$
here $\delta_{i,j}$ is the Kronecker symbol.
\end{lemma}
\noindent {\em Proof.} We let
$$\tilde{v}_i=\prod_{j=1,...,3p-1; j \neq i, j \neq
2p-i}([\omega]-h_{j,1}),$$ which is certainly non-zero (since it
has a nontrivial action on an irreducible $A(\trip)$-module). Fix
an index $i$, $1 \leq i \leq p-1$. We have to show that
$\lambda_i$ and $\nu_i$ exist and that $\lambda_i h_i +\nu_i \neq
0$, which is sufficient to argue that $v_i \neq 0$. From the very
definition it is easy to see that for every $f(x) \in
\mathbb{C}[x]$ we have
$$f([\omega])*w_i=f(h_i)w_i$$ and clearly
$$([\omega]-h_{i,1})* \tilde{v}_i=w_i.$$
The last two formulas yield \bea && f([\omega]) *
\tilde{v}_i=(f(h_i)+([\omega]-h_{i,1})r([\omega]))* \tilde{v}_i
\nonumber \\
&& =f(h_{i,1})*\tilde{v}_i+r(h_{i,1})w_i, \nonumber \eea where
$r(x)$ is the unique polynomial satisfying
$$f(x)=r(x)(x-h_{i,1})+f(h_{i,1}).$$
Specialize now
$$f(x)=\prod_{j=1,...,3p-1; j \neq i, j \neq
2p-i}(x-h_{j,1}),$$ so that
$$\tilde{v}_i *
\tilde{v}_i=f(h_{i,1})\tilde{v}_{i}+r(h_{i,1})w_i,$$ \be
\label{wi} \tilde{v}_i
* w_i =f(h_{i,1})w_i, \ee and consequently
$$v_i * v_i =(\lambda_i[\omega]+\nu_i)^2 (\tilde{v}_i *
\tilde{v}_i)$$
$$=f(h_{i,1})((\lambda_i h_{i,1}+\nu_i)^2)\tilde{v}_{i}+2
\lambda_i f(h_{i,1})(\lambda_i h_{i,1}+\nu_i)
w_i+r(h_{i,1})(\lambda_i h_{i,1}+\nu_i)^2w_i$$
$$=f(h_{i,1})((\lambda_i h_{i,1}+\nu_i)^2)\tilde{v}_{i}+
(\lambda_i h_{i,1}+\nu_i)(2 \lambda_i
f(h_{i,1})+r(h_{i,1})(\lambda_i h_{i,1}+\nu))w_i$$
$$=(K+L([\omega]-h_{i,1}))*\tilde{v}_i$$
where we let
$$K=f(h_{i,1})((\lambda_i h_{i,1}+\nu_i)^2)$$
and
$$L=(\lambda_i h_{i,1}+\nu_i)(2 \lambda_i f(h_{i,1})+r(h_{i,1})(\lambda_i h_{i,1}+\nu_i)).$$
Since we want
$$v_i * v_i =v_i,$$ by comparing the coefficients we get a system
$$L=\lambda_i, \ \ \ K-L h_{i,1}=\nu_i.$$
From the second equation we obtain \be \label{lmu} \lambda_i
h_{i,1}+\nu_i=\frac{1}{f(h_{i,1})} \neq 0. \ee From this relation,
after some computation, we get
$$\lambda_i=-\frac{r(h_{i,1})}{f(h_{i,1})^2},$$
and
$$\nu_i=\frac{1}{f(h_{i,1})}+\frac{h_{i,1} r(h_{i,1})}{f(h_{i,1})^2}.$$
From (\ref{wi}) and (\ref{lmu}) we clearly have
$$v_i * w_i=w_i * v_i=(\lambda_i h_{i,1}+\nu_i)(\tilde{v}_i * w_i)=w_i.$$
\epfv

\begin{example}
For instance, for $p=2$ and $i=1$, we have $h_{1,1}=0$ and
$$\tilde{v}_1=([\omega]+1/8)([\omega]-3/8)([\omega]-1).$$
From the previous lemma we have $r(0)=13/64$, $f(0)=3/64$
$$\lambda_1=-\frac{832}{9}, \ \  \nu_i=\frac{64}{3},$$
and hence
$$v_1=-\frac{64}{9}(13[\omega]-3)([\omega]+1/8)([\omega]-3/8)([\omega]-1).$$
\end{example}

\begin{remark} {\em It is not at all clear that $w_i \neq 0$!}.
\end{remark}

For $2p \le i \le 3p-1$, let $\mathbb{M}_{h_{i,1}}$ be the vector
space spanned by the vectors
$$A^{(i)}= C_p \prod_{j=2p, \dots, 3p-1, j\ne i} ([\omega]-h_{j,1}) * P([\omega])  \
, $$
$$B^{(i)}= \prod_ {j=2p, \dots, 3p-1, j\ne i} ([\omega]-h_{j,1}) * \
[H], $$
$$C^{(i)}= \prod_{j=2p, \dots, 3p-1, j\ne i} ([\omega]-h_{j,1})  *\
[E], $$
$$D^{(i)}= \prod_{j=2p, \dots, 3p-1, j\ne i} ([\omega]-h_{j,1})  *\
[F]. $$
It is easy to see that $\mathbb{M}_{h_{i,1}}$ is a nontrivial vector
space. From Theorem \ref{zhu-relations} follow that $[\omega]$ acts
on $\mathbb{M}_{h_{i,1}}$ by the scalar $h_{i,1}$.

\begin{lemma} \label{to-be-written}
\item[(i)] For every $i \in {\N}$, $2p \le i \le 3p-1$,
$\mathbb{M}_{h_{i,1}}$ is an ideal in $A(\trip)$ isomorphic to the
matrix algebra $M_2(\mathbb{C})$.

\item[(ii)] $\bigoplus_{i=2p} ^{3p-1} \mathbb{M}_{h_{i,1}}$ is an
ideal in $A(\trip)$.
\end{lemma}
\noindent {\em Proof.} Theorem \ref{zhu-relations} implies that
$\mathbb{M}_{h_{i,1}}$ is an ideal in $A(\trip)$. Moreover,
irreducible $A(\trip)$--module $\Pi(3p-i) (0)$ is also an
irreducible $2$--dimensional $\mathbb{M}_{h_{i,1}}$--module. Thus
we have a homomorphism from  $ \mathbb{M}_{h_{i,1}}$ to
$M_2(\mathbb{C})$. But this has to be an isomorphism, because the
module is irreducible.  This proves (i). Assertion (ii) follows
from (i) and the fact that $[\omega]$ acts on
$\mathbb{M}_{h_{i,1}}$ by the scalar $h_{i,1}$. \epfv

\vskip 5mm

Define :
$$ v_p = \prod_{j=1, \dots,3p-1, j \ne p}([\omega]-h_{j,1}); \quad  \C_{h_{p,1}}=\mathbb{C}v_{p}. $$
We have:
$$v_p * v_i=v_p  *w_i= 0, \ \ i=1,...,p-1$$
and
$$v_p * [\omega]=h_{p,1}
v_p.$$ Also, it is easy to see (by using Theorem
\ref{zhu-relations},(v)) that
$$X * v_p=v_p * X=0, \ \ X \in \{A^{(i)},B^{(i)},C^{(i)},D^{(i)} \}.$$
Thus, we have
\begin{lemma}
$ \mathbb{C}_{h_{p,1}}$ is a one-dimensional ideal in $A(\trip)$.
\end{lemma}

Let $\mathbb{I}_{h_{i,1}}$ be the ideal in $A(\trip)$ spanned by
$v_i$ and $w_i$.  We expect the following to be true.
\begin{conjecture} \label{zhu-conj}
Each $\mathbb{I}_{h_{i,1}}$ is a two-dimensional ideal.
%which
%implies
%$${\rm dim}(A(\trip))=6p-1.$$
\end{conjecture}

%The previous theorem is unsatisfactory because we do not know
%whether the ideal $I_{h_{i,1}}$ is two-dimensional or only
%one-dimensional. We of course expect a stronger result
%\begin{conjecture} \label{zhu-conj}
%Each $\mathbb{I}_{h_{i,1}}$ is a two-dimensional ideal.
%%which
%implies
%$${\rm dim}(A(\trip))=6p-1.$$
%\end{conjecture}
We will return to Conjecture \ref{zhu-conj} in the next section.

Now, we summarize the results from this section.
\begin{theorem} \label{pre-zhu-desc}
%Assume that Conjecture \ref{zhu-conj} holds.
  Zhu's algebra $A(\trip)$ decomposes as
a direct of sum of ideals
$$A(\trip)=\bigoplus_{i=2p}^{3p-1}  \mathbb{M}_{h_{i,1}} \oplus \bigoplus_{i=1}^{p-1} \mathbb{I}_{h_{i,1}} \oplus \mathbb{C}_{h_{p,1}},$$
where $\mathbb{M}_{h_{i,1}}$, $\mathbb{I}_{h_{i,1}}$ and
$\mathbb{C}_{h_{p,1}}$ are as above.
% and $\mathbb{I}_{h_{i,1}}$ is
%the ideal spanned by $v_i$ and $w_i$.
%
 Assume that Conjecture
\ref{zhu-conj} holds. Then $${\rm dim}(A(\trip))=6p-1.$$
% ideal with
%$$v_{-(p-1)^2/(4p)}=\prod_{j=1}([\omega]-h_{j,1}).$$
\end{theorem}
%\begin{proposition}
%There exists a nonzero polynomial $q(x)$ of degree $\leq p$ such
%that
%$$\prod_{i=2p}^{3p-1} ([w]-h_{i+1,1})[X]=0, \ \ X \in \{ E, F, H \}$$
%for
%$$q(\bar{\omega})H=q(\bar{\omega})E=q(\bar{\omega})F=0,$$
%which holds in $A(\trip)$. Consequently, we know that
%$$H_{-2}F=\bar{H} \bar{F}=$$
%\end{proposition}

%\noindent {\em Proof.} Firstly, we observe that $v=F_{-2}H$ is in
%the Virasoro submodule generated by $F$, so that there exists $a
%\in U(Vir)_{<0}$ such that $F_{-2}H=a \cdot H$, so there is $q \in
%\C[x]$ such that $q(\bar{\omega})H=0$. We will derive some
%properties of $q(\bar{\omega})$ first. Since ${\rm deg}(a)=2p$, it
%is not hard to see that
%$$\bar{a H}=\lambda \bar{\omega}^p \bar{H},$$ for some $\lambda
%\in \mathbb{C}$. The hard part is to show that $\lambda \neq 0$.
%In order to see that we have to go back to $A(\trip)$. We know
%that
%$$[F_{-2}H]=\tilde{q}([\omega])H,$$
%where $\tilde{q}(x)$ a polynomial of degree less or equal to $p$.
%If this polynomial is of deg

%Also, $\bar{v}=\bar{F_{-2}H}=0$.

%\epfv

\section{Logarithmic $\mathcal{W}(p)$-modules}

In this section we prove the existence of certain logarithmic
$\trip$-modules needed for the description of ideals
$\mathbb{I}_{h_{i,1}}$.

Let us recall that a logarithmic module for a vertex operator
algebra is a weak $\trip$-module which admits a decomposition into
generalized $L(0)$-subspaces. Since $\trip$ satisfies the
$C_2$-property a result of Miyamoto \cite{Miy} implies that every
weak $\trip$-module is logarithmic. Nontrivial logarithmic modules
(in the sense that they admit nontrivial Jordan blocks) are always
reducible. It is a priori not clear that non-trivial logarithmic
$\trip$-modules actually exist, so it is still an open problem to
construct them explicitly (We feel that the approach from
\cite{FFHST} and \cite{M1} might be useful for those purposes).
For $p=2$ a single logarithmic module can be constructed
explicitly by using {symplectic fermions} as shown in \cite{GK1}
(cf. also \cite{Abe}). An additional difficulty with logarithmic
modules is that they might involve Jordan blocks of larger size
deeper in the grading compared with those on the top. Thus, the
lowest weight subspace, being a $A(\trip)$-module, does not carry
enough information about the module itself and it cannot be used
to give even an upper bound on the size of Jordan blocks. Luckily
there are higher analogs of Zhu's algebra, denoted by $A_n(V)$, $n
\geq 1$ that control the representation theory below the top
component. We plan to return to $A_n(\trip)$ in forthcoming
publications \cite{AdM2}.

Even though in this paper we did not develop proper algebraic
tools to study logarithmic modules, we do have analytic tools
stemming from generalized graded traces of weak $\trip$-modules
\cite{Miy}. So in what follows we shall try to prove the existence
of logarithmic modules by using an indirect approach. Our proof is
in the spirit of the proof of existence of $g$-twisted sectors for
$g$-rational vertex operator algebras obtained by Dong, Li and
Mason in \cite{DLM}. In their approach, non-triviality of
$g$-twisted sectors (i.e., existence of a $g$-twisted module) is
proven by using modular invariance. Similarly, here we employ
Miyamoto's modular invariance of {\em pseudotraces} for vertex
algebras satisfying $C_2$-cofiniteness condition. We should say
here that Miyamoto's result provides us only with {\em some}
logarithmic modules. The hard part is to construct {\em enough}
logarithmic $\trip$-modules with two-dimensional lowest weight
subspaces, such that Conjecture \ref{zhu-conj} holds true. We will
eventually show that this could be done, without too much effort,
for $p$ being a prime integer.

%\begin{proposition} If $\mathcal{W}(p)$ admits $p$ (nonisomorphic)
%logarithmic modules with two-dimensional top weight component of
%generalized conformal weight $h_{1,p}$, $i=1,...,p-1$, then the
%Conjecture \ref{zhu-conj} holds true.
%\end{proposition}
%As we already mentioned, in \cite{Miy} a logarithmic version of
%Zhu's main theorem \cite{Zhu} was obtained for vertex operator
%algebras satisfying the $C_2$-cofiniteness condition.

We first recall some notation and results from \cite{Miy}. For
every $n \geq 0$ the $n$-th Zhu's associative algebra is defined
as $A_n(V)=V/O_n(V)$, where $O_n(V)$ is spanned by
$${\rm Res}_x \frac{(1+x)^{n+{\rm deg}(u)}}{x^{2n+2}} Y(u,x)v,$$
where $u, v \in V$ (with $u$ homogeneous), and the $n$-th product
$*_n$ in $A_n(V)$ is defined similarly as for $A_0(V)=A(V)$ (see
\cite{Miy} for details). If $V$ is $C_2$-cofinite then all
$A_n(V)$ are in fact finite-dimensional associative algebras. If
$W=\coprod_{n \geq 0} W_n$ is an $\mathbb{Z}_{\geq 0}$-graded weak
$V$-module than $W([0,n])=\oplus_{ 0 \leq i \leq n} W_n$ is a
$A_n(V)$-module and $o([a] *_n [b])=o(a)o(b)$ holds on $W([0,n])$.
Suppose that $T$ is an $A_n(V)$-module. Then one forms a
generalized Verma $V$-module $W_T(n)$ (a weak $V$-module), whose
$n$-th graded piece $W(n)$ is $T$. Another important gadget is a
{\em pseudotrace} ${\rm tr}_{W_T(n)}^{\phi}$, where $\phi$ is a
symmetric map interlocked with $W_{T}(n)$ (for definitions see
\cite{Miy}). The main result of Miyamoto is then
\begin{theorem} \label{miy-mod} Suppose that $V$ is $C_2$-cofinite VOA.
Let $n$ be large enough and $W^1$,...,$W^m$ be $n$-th generalized
Verma $V$-module, interlocked with symmetric functions $\phi_i$.
Then the vector space spanned by ${\rm tr}^{\phi_i}_{W^i}
q^{L(0)-c/24}$, is modular invariant. In addition, in the $\tau$
expansion of ${\rm tr}^{\phi_i}_{W^i} q^{L(0)-c/24}$ all
coefficients are ordinary (pseudo)characters.
\end{theorem}

Since $L(0)$ does not act semisimply in general in all these
formulas with pseudotraces we use
$$q^{L(0)}=\sum_{k \geq 0} (2 \pi i \tau)^k \frac{L_{n}^k(0)}{k!}
q^{L_{ss}(0)},$$ $L(0)=L_{ss}(0)+L_n(0)$, where $L_n(0)$ is the
nilpotent and $L_{ss}(0)$ the semismple part of $L(0)$.

Now, we take $V$ to be $\trip$ and discuss first the irreducible
$\trip$-characters. We denote by
$$\eta(q)=q^{1/24} \prod_{n \geq 1}(1-q^n),$$
the Dedekind $\eta$-function and by
\bea \label{theta-const} \theta_{i,p}(q)&=&\sum_{n \in \mathbb{Z}}
q^{\frac{(2pn+i)^2}{4p}},
\\ (\partial \theta)_{i,p}(q)&=& \sum_{n \in \mathbb{Z}} (2pn+i)
q^{\frac{(2pn+i)^2}{4p}}, \eea certain theta constants and their
derivatives (these are in fact modular forms of weight $1/2$ and
$3/2$ for some congruence subgroups, respectively). From Theorem
\ref{str-ff-1} and Theorem \ref{str-ff-2} it is not hard to see
(by using well-known formulas for the characters of irreducible
Virasoro modules) that in fact (cf. \cite{FHST}, \cite{F1},
\cite{FG}) \bea \label{irr-char} {\rm tr}_{\Lambda(i)}
q^{L(0)-c_{p,1}/24} & = & \frac{1}{p \eta(q)}\left( i
\theta_{p,p-i}(q)+ (\partial
\theta_{p,p-i})(q)\right), \\
{\rm tr}_{\Pi(i)} q^{L(0)-c_{p,1}/24} & = & \frac{1}{p
\eta(q)}\left( i  \theta_{p,i}(q)- (\partial
\theta_{p,i})(q)\right). \eea

In what follows, we will say that two $q$-series $f \in q^a
\mathbb{C}[[q]]$ and $g \in q^b \mathbb{C}[[q]]$ have {\em
compatible} $q$-expansions if $a-b \in \mathbb{Z}$.

\begin{lemma} \label{lowest} Let $p \geq 2$ be a prime integer.
Then for $i \neq j$, $i,j \in \{1 \leq i \leq p \} \cup \{ 2p \leq
i \leq 3p-1 \}$,
$$h_{i,1}-h_{j,1} \in \mathbb{Z}$$
if and only if $i-j=2p$, $1 \leq i \leq p-1$ (or $i-j=-2p$).
\end{lemma}
\noindent {\em Proof.} The statement clearly holds for $p=2$, so
we may assume $p \geq 3$. Since $p$ is prime
$$h_{i,1}-h_{j,1}=\frac{(i-p)^2-(j-p)^2}{4p} \in \mathbb{Z},$$
implies that $p$ divides $i+j-2p$ or $i-j$. In the former case the
only possibilities are $i+j-2p=kp$, $k \in \{-1,0,1,2,3 \}$. If
$i+j-2p=\pm p$, then $i-j=4l$, which implies $2i=4l \pm 3p$,
having no solution. Similarly, if $i+j-2p=3p$. The case $i+j=2p$
and $i+j=4p$ have no solution for $i \neq j$. Therefore $p$
divides $i-j$, where the only possibility is $i-j=2p$ (or
$i-j=-2p$) or $i-j=p$. If $i-j=\pm p$, then $i+j-2p=p$, which
reduces to a previous case. Thus, we are left with $i-j= \pm 2p$.
The proof follows. \epfv

An important consequence of the previous lemma is that, for $p$
prime, the only irreducible characters which have compatible
$q$-expansions are $\Lambda(i)$, with the lowest conformal weight
$h_{i,1}$ and $\Pi(p-i)$, with the lowest conformal weight
$h_{2p+i,1}$, $i=1,...,p-1$. Notice also that \be \label{less}
h_{i,1} < h_{2p+i,1}, \ \ 1 \leq i \leq p-1. \ee

Next result is a consequence of modular properties of irreducible
$\trip$ characters (cf. \cite{F1} for instance). Its proof is
well-recorded so we omit the proof here (see \cite{F1} for
instance).
\begin{lemma} \label{flohr-95} The vector space spanned by $2p$-irreducible $\trip$-characters
is $2p$-dimensional. The $SL(2,\mathbb{Z})$ transforms of these
characters closes $3p-1$-dimensional $SL(2,\mathbb{Z})$-module
with a basis formed by irreducible characters ${\rm
tr}_{\Lambda(i)} q^{L(0)-c/24}$, ${\rm tr}_{\Pi(i)}
q^{L(0)-c/24}$, $1 \leq i \leq p$ and
$$2 \pi i \tau \frac{(\partial \theta)_{p,i}(q)}{\eta(q)}, \ i=1,2,...,p-1.$$
\end{lemma}

Now we are ready to prove the main theorem in this section.
\begin{theorem} \label{some-logar} For every prime $p$ and $i \in \{1,...,p-1 \}$, $\trip$ admits a
logarithmic module with a two-dimensional lowest weight subspace
of generalized conformal weight $h_{i,1}$.
\end{theorem}
\noindent {\em Proof.} Since the triplet is $C_2$-cofinite, every
module is $\mathbb{Z}_{\geq 0}$-gradable and logarithmic (which
includes also ordinary modules). Let $W$ be a weak $\trip$-module
with a top component $W(0)$. Because of $f_p([\omega])=0$, which
holds in Zhu's algebra we have $W(0)=\oplus_{i=1}^{3p-1}
W_{h_{i,1}}(0)$, where $W_{h_{i,1}}(0)$, is $h_{i,1}$-primary
component of $W(0)$, that is $W_{h_{i,1}}$ is annihilated by
$(L(0)-h_{i,1})^2$.

From Lemma \ref{flohr-95} we know that \be \label{char-tau0} 2 \pi
i \tau \frac{(\partial \theta)_{p,i}(q)}{\eta(q)} \in \tau
q^{h_{i,1}-c_{p,1}/24}(b_i+ q\mathbb{C}[[q]]), \ \ b_i \neq 0 \ee
is an $SL(2,\mathbb{Z})$-transform of the set of irreducible
characters. Fix an index $i$. According to the discussion
preceding Theorem \ref{miy-mod}, we have \be \label{char-tau1} 2
\pi i \tau \frac{(\partial \theta)_{p,i}(q)}{\eta(q)}=\sum_{j=1}^k
{\rm tr}^{\phi_j}_{W^j} q^{L(0)-c_{p,1}/24}, \ee for some
pseudotraces ${\rm tr}^{\phi_j}_{W^j}$, depending on $i$, with
symmetric maps $\phi_j$ evaluated on weak modules $W^j$ (if we
multiply a pseudotrace with a constant we obtain another
pseudotrace). We may assume that on the right hand-side in
(\ref{char-tau1}) we only have pseudotraces with the coefficients
in the $\tau$-expansion being $q$-series compatible with the left
hand-side in (\ref{char-tau1}), so let
$${\rm tr}_{W^j}^{\phi_j} q^{L(0)-c_{p,1}/24} \in
\sum_{m=0}^r \tau^m q^{\tilde{h}-c_{p,1}/24} \mathbb{C}[[q]],$$
for some $\tilde{h}$, with $\tilde{h}-{h_{i,1}} \in \mathbb{Z}$.
According to Lemma \ref{lowest} we may assume
$\tilde{h}={h_{i,1}}$, because there will be no $\trip$-module
with generalized conformal weight less than $h_{i,1}$. Of course,
because of (\ref{less}), characters of modules with lowest
conformal weight $h_{2p+i,1}$ are already contained in
$q^{h_{i,1}}\mathbb{C}[[q]]$. Thus we may assume that all $W^j$
have the lowest generalized conformal weights $h_{i,1}$ or
$h_{2p+i,1}$. We claim that for at least one $W^j$, the lowest
weight subspace $W^j(0)$ admits at least one nontrivial
$L(0)$-Jordan block of length two of generalized conformal weight
$h_{i,1}$. As we already mentioned the property
$$(L(0)-h_{i,1})^{2}W^j(0)=0, \ i \geq 0,$$
guarantees that $L(0)$ Jordan blocks cannot be of size strictly
bigger than $2$ (but this does not rule out existence of Jordan
blocks of larger size below the top!). Suppose that there is no
such $W^j$, that is suppose that the lowest weight subspace
$W^j(0)$ is $L(0)$-diagonalizable or $W^j(0)$ is trivial. Then
from \cite{Miy}, for every $j$ we have \bea \label{miy-trick} {\rm
tr}_{W^j}^{\phi_j} q^{L(0)-c_{p,1}/24}&=&{\rm tr}^{\phi_j}_{W^j}
\sum_{k \geq 0} \frac{(2 \pi i \tau)^k
(L(0)-L_{ss}(0))^k}{k!} q^{L_{ss}(0)-c_{p,1}/24} \\
&=&\sum_{k \geq 0} (2 \pi i \tau)^k {\rm
tr}^{\tilde{\phi}_k}_{W^j/R_k} q^{L_{ss}(0)-c_{p,1}/24}, \eea
where $R_0=0$ and $R_k$, $k \geq 0$ are certain submodules of
$W_j$. But if $W^j(0)$ is $L(0)$-diagonalizable then there is no
$\tau$-term in (\ref{miy-trick}), because $(L(0)-L_{ss}(0))$
annihilate all of $W_i(0)$, so in the $(q,\tau)$-expansion of
${\rm tr}_{W^j}^{\phi_j} q^{L(0)-c_{p,1}/24}$ there will be no
term of the form $\tau q^{h_{i,1}}$, contradicting to
(\ref{char-tau0}) and (\ref{char-tau1}). \epfv

Here is an important consequence of the previous theorem.
\begin{theorem} \label{zhu-p}
The Conjecture \ref{zhu-conj} holds for every prime $p$.
\end{theorem}
\noindent {\em Proof.} In view of Theorem \ref{pre-zhu-desc} and
Conjecture \ref{zhu-conj} we only have to show that each
$\mathbb{I}_{h_{i,1}}$ is two-dimensional. Since each $v_i$ is an
idempotent and $w_i^2=0$, $v_i$ and $w_i$ are not proportional.
Thus, it is sufficient to show $w_i \neq 0$, or that $w_i$ acts
non-trivially on a $A(\trip)$-module. Since $\trip$ admits $p-1$
logarithmic modules, let's call them $\Lambda_2(i)$, $i=1,..,p-1$,
such that $w_i \in A(\trip)$ acts nontrivially on the top
component $\Lambda_2(i)(0)$ the proof follows. \epfv

The $p=2$ case of Theorem \ref{zhu-p} has been verified by Abe
\cite{Abe} by using explicit fermionic construction of logarithmic
modules, obtained previously by physicists.

In general, we are able to prove the following result.

\begin{proposition} \label{gen-log-p}
For every $p \geq 2$, the triplet $\trip$ admits a logarithmic
module of lowest conformal weight $h_{p-1,1}$. Moreover,
$\mathbb{I}_{h_{p-1,1}}$ is two-dimensional.
\end{proposition}
\noindent {\em Proof.} The main argument is similar to the one
used in the proof of Theorem \ref{some-logar} so we omit details.
First observe that $h_{p,1}=\frac{-(p-1)^2}{4p}$ is the smallest
conformal weight among all irreducible $\trip$-modules. Since
$h_{p-1,1}=\frac{1-(p-1)^2}{4p}$ is not congruent to $h_{p,1}$ mod
$\mathbb{Z}$ for any $p \geq 2$, and $h_{p-1,1}<h_{i,1}$ for other
$i$, then as in Theorem \ref{some-logar} it follows that $\trip$
admits a logarithmic module with a two-dimensional generalized
lowest conformal weight $h_{p-1,1}$. Consequently,
$\mathbb{I}_{h_{p-1,1}}$ is two-dimensional. \epfv

We finish with an expected conjecture
\begin{conjecture} There are no logarithmic $\trip$-modules
admitting $L(0)$ Jordan blocks of size three or more.
\end{conjecture}

\section{Final remarks and future work}

In this section we gather some problems and open question that we
shall address in our future publications.

\begin{itemize}

\item[(i)] The problem of constructing logarithmic $\trip$-modules
tops the list of our future directions. It seems to us that this
problem hasn't been solved even in the physics literature, except
for the $p=2$ case. We feel that approach in \cite{FFHST} and
\cite{M1} might be useful for these purposes.
%A closely related
%problem is decomposition of the category of $\trip$-modules into
%{\em linkage} classes (cf. \cite{FHST}).

\item[(ii)] Vertex operator algebra $\trip$ can be considered as a
$\goth{sl}_2 \times {\rm Vir}$--module (cf. \cite{F1},\cite{F3},
\cite{FGST1}, \cite{FGST2}, \cite{FGST3}). This also follows from
our vertex-algebraic approach.
Define the following operators acting on $\trip$:
$$ e= Q, \quad h = \frac{\a(0)}{p}.$$
 From Theorem
\ref{str-ff-1} and Proposition \ref{stgen} follow that  there is a
unique operator $f\in \mbox{End}(\trip)$ which commutes with the
action of the Virasoro algebra such that
$$ f  e^{-n\a} =  0, \quad f Q^{j} e^{-n\a} =  -j (j -1-2n ) Q^{j-1}
e^{-n\a}, \ 1\le j \le 2n. $$
Therefore  $\trip$ is an $\goth{sl}_2 \times L(c_{p,1},0)$--module
and
$$\trip = \bigoplus_{n= 0} ^{\infty} W_{2n+1} \otimes
L(c_{p,1},n^2 p + np -n)$$
where $W_{2n+1}$ is a $(2n+1)$--dimensional $\goth{sl}_2$-module.
This implies that the Lie group $PSL(2, {\C})$ acts on the vertex
operator algebra $\trip$ as an automorphism group.

 Let
$\Gamma \subset PSL(2,\C)$ be any finite group. Then one can
consider the fixed point subalgebra $\trip ^{\Gamma}$. We think
that it is important and interesting problem to investigate
$C_2$--cofiniteness  of these subalgebras of $\trip$.

\item[(iii)] Finally, there is a fair amount of work needed to
determine the fusion rules of logarithmic and non-logarithmic
$\trip$-modules in the context of vertex algebras. Since the
category of $\trip$-modules has a natural braided tensor category
structure \cite{CF}, \cite{HLZ1}-\cite{HLZ2} it is natural to look
for an already existing model for this category. Physicists have
provided a beautiful conjecture in this direction: {\em the tensor
category of $\trip$-modules is equivalent to the tensor category of
$\overline{{U}_q(\goth{sl}_2)}$-modules, $q=e^{\frac{\pi i}{p}}$}.
This has been verified for $p=2$ in \cite{FGST1}, \cite{FGST3}.

\end{itemize}

\section{Appendix}

\subsection{Proof of Theorem \ref{zhu-relations}, (iii)}

From the Vir-module structure of $\trip$ it follows that there
exists $a \in \mathcal{U}(Vir_{\leq -1})$ such that
$$H * F =a.F.$$
Since $$[H*F]=[H]*[F]$$ in $A(\trip)$ and
%
%${\rm deg}(H*F)=4p-2$, (promijenjeno)
%
$${\wt}(H_{-1}F) = 4p-2$$
%
%maximal weight of homogeneous elements which appear in the
%expression for  $H*F$ is $4p-2$,
 we conclude that  there exists a
polynomial $q(x) \neq 0$, ${\rm deg}(q(x)) \leq p-1$, such that
$$[H]*[F]=-q([\omega])*[F].$$
Also, we know that $\Pi(i)$ is an $A(\trip)$-module with
$2$-dimensional top weight subspaces $\Pi(i)(0)$ of lowest
conformal weights $h_{i,1}$, $i=2p,...,3p-1$. Thus, we obtain the
relation \be \label{top-rel}
[H]*[F]|_{\Pi(i)(0)}=-q([\omega])*[F]|_{\Pi(i)(0)}. \ee It is not
hard to see that this relation uniquely determines $q(x)$ and that
${\rm deg}(q(x))=p-1$ (this will be proven below). As in Section
3, it is convenient to switch to "charge variable" $t$. Thus, we
define \be \label{hp}
 H_p(t)=-q \left(\frac{t(t-2p+2)}{4p}\right)
\in \mathbb{C}[t]. \ee Now, $[H]$ acts on $M(1,t)(0)$ as
multiplication with ${t \choose 2p-1}$, and $[\omega]$ acts as
multiplication with $\frac{t(t-2p+2)}{4p}$. By using Lagrange
interpolation theorem and relation (\ref{top-rel}) and (\ref{hp})
it is not hard to prove the following result:
\begin{lemma}
We have
$$ H_p(t)=\left( \frac{1}{(2p-1)!}\prod_{i=2p-1}^{3p-2} (t-i)(t+i-2p+2)
\right) \cdot $$ $$\cdot \left(\sum_{i=2p-1}^{3p-2}
\frac{(-1)^{p-i}(i!)^2}{(i-2p+1)!^2(p+i)!(3p-i-2)!}
\left(\frac{1}{t-i}-\frac{1}{t+i-2p+2}\right)\right).$$
\end{lemma}

Here are the first few $H_p(t)$ polynomials:

\bea H_2(t)&=&\frac{3}{5}t^2-\frac{6}{5}t-\frac{4}{5} \nonumber \\
H_3(t)&=&\frac{5}{84}t^4-\frac{10}{21}t^3+\frac{55}{84}t^2+\frac{25}{21}t+1
\nonumber \\
H_4(t)&=&
\frac{35}{15444}t^6-\frac{35}{858}t^5+\frac{3395}{15444}t^4-\frac{245}{1287}t^3-\frac{1897}{3861}t^2
-\frac{3136}{1287}t-\frac{200}{143} \nonumber \eea and the
corresponding $q(x)$ polynomials are then \bea &&
q(x)=-(\frac{24}{5}x- \frac{4}{5}), \ \ {\rm for} \ p=2, \nonumber \\
&& q(x)=-(\frac{60}{7}x^2-\frac{25}{7}x+1), \ \ {\rm for} \ p=3, \nonumber \\
&&
q(x)=-(\frac{35840}{3861}x^3-\frac{2240}{351}x^2+\frac{25088}{3861}-\frac{200}{143}),
\ \ {\rm for} \ p=4. \nonumber \eea

As far as we can tell $q(x)$ polynomials do not admit nice
factorization so we need a different approach to show that
$q(h_{i,1}) \neq 0$. Let
$$H_p(t)=\frac{Pr_p(t)}{(2p-1)!} (S_p(t) + \tilde{S}_p(t)),$$
where
$$\tilde{S}_p(t)=\sum_{i=2p-1}^{3p-2}
\frac{-(-1)^{p-i}i!^2}{(t+i-2p+2)(i-2p+1)!^2(p+i)!(3p-i-2)!}$$
$${S}_p(t)=\sum_{i=2p-1}^{3p-2}
\frac{(-1)^{p-i} i!^2}{(t-i)(i-2p+1)!^2(p+i)!(3p-i-2)!}$$ and
$$Pr_p(t)=\prod_{i=2p-1}^{3p-2}
(t-i)(t+i-2p+2).$$
\begin{proposition} For every $p \geq 2$,
\be \label{form} H_p(t) \neq 0, \ \ t \in [0,2p-2] \cap
\mathbb{Z}. \ee
\end{proposition}
\noindent {\em Proof:} Notice that $Pr_p(t) \neq 0$ for $t \in
\{0,1,...,2p-2\}$, so we only have to consider
$$A_p(t)=S_p(t) + \tilde{S}_p(t).$$
It is easy to see that
$$A_p(2p-2-t)=A_p(t)$$
for an arbitrary value of parameter $t$. Thus, relation
(\ref{form}) will follow once we prove: $A_p(t) \neq 0$ for $t \in
[0,p] \cap \mathbb{Z}$.

It turns out that there isn't nice closed formula $A_p(t)$, but
$A_p(0)$ and $A_p(1)$ can be computed explicitly. By using
standard hypergeometric summation techniques we obtain:
$$A_p(0)=-1/2\frac{4^p (2p-1){p-3/2 \choose p}^2}{{2p-3/2 \choose
p}},$$
$$A_p(1)=-\frac{4^{p-1}(2p-1)(2p^2-1){p-3/2 \choose p}^2}{(p^2-1){2p-3/2 \choose
p}}.$$ Both $A_p(0)$ and $A_p(1)$ are evidently less than zero.

Furthermore, by using straightforward computation we also have a
degree two recursion
$$A_p(t)=\frac{
(t-1)^2(3p-t)A_p(t-2)+2(t-p)(t^2-2pt+2p-2p^2)A_p(t-1)}{(t+p)(t+1-2p)^2}.$$
Now it is not hard to see that for $t \in [0,p]$,
$$(t-p)(t^2-2pt+2p-2p^2) \geq 0$$
and of course
$$(t-1)^2(3p-t) \geq 0.$$
Now, inductively it easily follows that $A_p(t)<0$ for $t \in
[0,2p-2] \cap \mathbb{Z}$. The proof follows. \epfv

From the previous proposition we get $q(h_{i,1}) \neq 0$.

%From
%$$[H]*[F]=q([w])*[F]$$ we obtain
%$$[H]^2*[F]=q([w])*[H]*[F],$$ and
%$$(C_p P([w])-p([w])^2)*[F]=0,$$
%where $C_p P([w])=[H]^2$ is as in \cite{A-2003}. Since
%$$p(h_{i,1})^2=C_p P(h_{i,1}), \ {\rm for} \ \ 2p \leq i \leq 3p-1,$$
%and
%$$q(h_{i,1}) \neq 0, \ {\rm for} \  \ i=1,..,2p-1,$$
%by part (iii), we have
%$$(P([w])-p([w])^2)*[F]=r([w]) \prod_{i=2p}^{3p-1}([w]-h_{i,1})*[F]=0,$$
%where $r([w])$ is a polynomial relatively prime with
%$$f_p([w])=\prod_{i=1}^{3p-1}([w]-h_{i,1}).$$ Since $f_p([w])$ is plainly zero in $A(W(p))$, then $r([w]) \neq 0$ is a
%unit in $A(W(p))$, so it can be cancelled from the left. \epfv
%Let
%$$h_{i,1}=\frac{i(i-2p+2)}{4p}$$
%(this is $h_{i+1,1}$ in your notes). Then, $q(h_{i,1}) \neq 0$,
%$i=0,..,2p-2$.


\begin{thebibliography}{EFHHNV}


\bibitem[1]{Abe} T. Abe,
A ${\Bbb Z}\sb 2$-orbifold model of the symplectic fermionic vertex
operator superalgebra. {\em Math. Z.} {\bf 255} (2007), 755--792.

\bibitem[2]{ABD} T. Abe, G. Buhl and C. Dong, Rationality, regularity, and $C_2$-cofiniteness,
{\em Trans. Amer. Math. Soc.} {\bf 356}, (2004), 3391-3402.

%Classification of irreducible modules for the vertex operator
%algebra $V\sp +\sb L$: general case. {\em J. Algebra} {\bf 273}
%(2004), 657--685.


\bibitem[3]{A1} D. Adamovi\' c, Representations of the vertex
algebra ${\W}_{ 1 + \infty}$ with a negative integer central
charge, Comm. Algebra 29 (2001), no. 7, 3153-3166.


\bibitem[4]{A-2003} D. Adamovi\'{c},
Classification of irreducible modules of certain subalgebras of
free boson vertex algebra, J. Algebra 270 (2003) 115-132.

\bibitem[5]{AdM} D. Adamovi\'c and A. Milas, Logarithmic intertwining operators
and $\mathcal{W}(2,2p-1)$-algebras, {\em Journal of Math. Physics}
{\bf 48}, 073503 (2007).

\bibitem[6]{AdM2} D. Adamovi\'c and A. Milas, in preparation.

%\bibitem[3]{Andrews} G.E. Andrews, $q$-Series: Their Development and Application in Analysis, Number Theory,
%Combinatorics, Physics, and Computer Algebra, CBMS Regional Conf.
%Ser. in Math., vol. 66, American Mathematical Society, Providence,
%RI, 1986.

\bibitem[7]{Ar} T. Arakawa, Representation Theory of W-Algebras,
{\em Inventiones Math.} {\bf 169} (2007), 219-320. {\tt
arxiv:math/0506056}.

%\bibitem[4]{DGK} V. Deodhar, O. Gabber and V.Kac, Structure of some categories of representations of infinite-dimensional Lie algebras.
%Adv. in Math. 45 (1982), no. 1, 92--116.

\bibitem[8]{CF} N. Carqueville and M. Flohr,  Nonmeromorphic operator product expansion
and $C\sb 2$-cofiniteness for a family of $\mathcal{W}$-algebras.
{\em J. Phys.} A {\bf 39} (2006), 951--966.

\bibitem[9]{dSK} A. De Sole and V. Kac, Finite vs. affine
$W$-algebras, {\em Japanese Journal of Math.} {\bf 1} (2006)
137-261; arxiv:math/0511055.


\bibitem[10]{D} C. Dong, Vertex algebras associated with even lattices,  J. Algebra
{\bf 160} (1993), 245-265.

%\bibitem [8]{DL} C. Dong and J. Lepowsky,
%Generalized vertex algebras and relative vertex operators,
%Birkh\"auser,  Boston, 1993.

\bibitem[11]{DLTYY} C. Dong, C.H. Lam, K. Tanabe, H. Yamada, K.
Yokoyama, $Z_3$ symmetry and $W_3$ algebra in lattice vertex
operator algebras, {\em Pacific Journal of Mathematics} {\bf 215}
(2004), 245-296.

\bibitem[12]{DLM} C. Dong, H. Li and G. Mason, Modular-invariance of trace functions in orbifold theory and generalized Moonshine.
{\em Comm. Math. Phys.} {\bf 214} (2000), 1--56.

\bibitem [13]{EFH} W. Eholzer, M. Flohr, A. Honecker, R.
H\"ubel, W. Nahm and R. Vernhagen, Representations of
${\W}$--algebras with two generators and new rational models,
Nucl. Phys. B 383 (1992), 249-288.

%\bibitem [FFr]{FFr} B. Feigin and E. Frenkel, Integrals of motion
%and quantum groups, in Integrable Systems and Quantum groups,
%Lecture Notes in Math. 1620, p.p. 349-418, Springer Verlag, 1996.

\bibitem [14]{FF} B. Feigin and D. B. Fuchs, Representations of the Virasoro algebra,
in {\em Representations of infinite-dimensional Lie groups and Lie
algebras}, Gordon andd Breach, New York (1989).

\bibitem [15]{FGST1} B.L. Feigin, A.M. Ga\u\i nutdinov, A. M. Semikhatov, and I. Yu Tipunin, I,
The Kazhdan-Lusztig correspondence for the representation category
of the triplet $W$-algebra in logorithmic conformal field theories.
(Russian) Teoret. Mat. Fiz. {\bf 148} (2006), no. 3, 398--427.

\bibitem [16]{FGST2} B.L. Feigin, A.M. Ga\u\i nutdinov, A. M. Semikhatov, and I. Yu Tipunin,
Logarithmic extensions of minimal models: characters and modular
transformations. {\em Nuclear Phys.} B {\bf 757} (2006), 303--343.

\bibitem [17]{FGST3} B.L. Feigin, A.M. Ga\u\i nutdinov, A. M. Semikhatov, and I. Yu Tipunin,
Modular group representations and fusion in logarithmic conformal
field theories and in the quantum group center. {\em Comm. Math.
Phys.} {\bf 265} (2006), 47--93.


\bibitem [18]{Fel} G. Felder, BRST aproach to minimal models,
{ Nucl. Phys.} B {\bf 317} (1989), 215-236.

\bibitem [19]{FFHST} J. Fjelstad, J. Fuchs, S. Hwang, A.M.
Semikhatov and I. Yu. Tipunin, Logarithmic conformal field
theories via logarithmic deformations, {  Nuclear Phys.} B {\bf
633} (2002), 379--413.


\bibitem[20]{F1} M. Flohr, On modular invariant partition functions of conformal field theories with logarithmic operators,
{ Internat. J. Modern Phys.} A {\bf 11 } (1996), 4147--4172.


\bibitem[21]{F3} M. Flohr, Bits and pieces in logarithmic conformal field theory. Proceedings of the School and Workshop on Logarithmic Conformal Field Theory and its Applications (Tehran, 2001),
{ Internat. J. Modern Phys.} A {\bf 18} (2003), 4497--4591.

\bibitem[22]{FG} M. Flohr and M. Gaberdiel, Logarithmic torus amplitudes. {\em J. Phys.} A {\bf 39} (2006),
1955--1967.

\bibitem[23]{Fr} E. Frenkel, Lectures on the Langlands Program and Conformal Field
Theory, arxiv:math/0512172.


\bibitem[24]{FrB} E. Frenkel and D. Ben-Zvi, {\em Vertex algebras and algebraic curves},
Mathematical Surveys and Monographs, 88, American Mathematical
Society, Providence, RI, 2001.

\bibitem[25]{FKRW} E. Frenkel, V. Kac, A. Radul and W.Wang, $W_{1+\infty}$ and $W(gl_N)$ with central charge $N$, { Comm.
Math. Phys.} {\bf 170} (1995), 337-357.

%\bibitem [17]{FHL}
%I. B. Frenkel, Y.-Z. Huang and J. Lepowsky, On axiomatic
%approaches to vertex operator algebras and modules, Mem. Amer.
%Math. Soc. {\bf 104}, 1993.

%\bibitem[18]{FZ}
%I. B. Frenkel and Y.  Zhu, Vertex operator algebras associated to
%representations of affine and Virasoro algebras,  Duke Math. J. 66
%(1992),  123--168.


\bibitem[26]{FHST} J. Fuchs, S. Hwang, A.M. Semikhatov and I. Yu. Tipunin,
Nonsemisimple Fusion Algebras and the Verlinde Formula, { Comm.
Math. Phys.} {\bf 247} (2004), no. 3, 713--742.

%\bibitem [FLM]{FLM}
%I. B. Frenkel, J. Lepowsky and A. Meurman,   Vertex Operator
%Algebras and the Monster,   Pure and Applied Math., Vol. {\bf
%134}, Academic Press, New York, 1988.


\bibitem[27]{Ga} M. Gaberdiel, An algebraic approach to logarithmic conformal field theory, Proceedings of the School and Workshop on Logarithmic Conformal Field Theory and its Applications (Tehran, 2001),
{ Internat. J. Modern Phys.} A {\bf 18} (2003), 4593--4638.

\bibitem[28]{GK1} M. Gaberdiel and H. G. Kausch, A rational logarithmic
conformal field theory, Phys. Lett B 386 (1996), 131-137,
hep-th/9606050

\bibitem[29]{GK2} M. Gaberdiel and H. G. Kausch, A local logarithmic
conformal field theory, Nucl. Phys. B 538 (1999), 631-658,
hep-th/9807091

\bibitem[30]{GR} M. Gaberdiel and I. Runkel, The logarithmic triplet theory with
boundary, {\em J.Phys.} A {\bf 39} (2006), 14745-14780, {\tt
hep-th/0608184}.

%\bibitem[24]{Gu} V. Gurarie, Logarithmic operators in conformal field theory, { Nuclear Phys.} B
%{\bf 410} (1993), 535--549.

\bibitem[31]{H} A. Honecker, Automorphisms of ${\W}$ algebras
and extended rational conformal field theories, { Nucl. Phys.} B
{\bf 400} (1993), 574-596.


\bibitem[32]{H1} Y.-Z. Huang, Vertex operator algebras, the Verlinde conjecture,
and modular tensor categories. {\em Proc. Natl. Acad. Sci.} USA {\bf
102} (2005), 5352--5356

\bibitem[33]{HLZ1} Y.-Z. Huang, J. Lepowsky and L. Zhang, A logarithmic generalization of tensor
product theory for modules for a vertex operator algebra, {\em
Internat. Journal of Math.} {\bf 17}, (2006) 975–-1012, {\tt
math.QA/0311235}.

\bibitem[34]{HLZ2} Y.-Z. Huang, J. Lepowsky and L. Zhang, Logarithmic
tensor product theory for generalized modules for a conformal vertex
algebra, Part I, {\tt math.QA/0609833}.

\bibitem[35]{Kac}  V. Kac, {\em Vertex algebras for beginners}, University Lectures Series,
Vol. 10, Providence, 1998.

\bibitem[36]{Ka1} H. G. Kausch, Extended conformal algebras
generated by multiplet of primary fields, { Phys. Lett.} {\bf 259}
B (1991), 448-455.

\bibitem[37]{Ka2} H. G. Kausch, Symplectic Fermions, { Nucl. Phys.} B {\bf
583} (2000), 513-541.

\bibitem[38]{KaW} H. G. Kausch and G. M. T. Watts, A study of
${\W}$--algebras by using Jacobi identities, { Nucl. Phys.} B {\bf
354} (1991), 740-768.

\bibitem[39]{LL} J. Lepowsky and H. Li, {\em Introduction to Vertex Operator
Algebras and Their Representations}, Progress in Mathematics, Vol.
227, Birkh\"auser, Boston, 2003.

\bibitem[40]{M0} A. Milas, Fusion rings for degenerate minimal models, { J. Algebra} {\bf 254} (2002), no. 2, 300--335.

\bibitem[41]{M2} A. Milas, Weak modules and logarithmic intertwining operators for vertex operator algebras. Recent developments in infinite-dimensional Lie algebras and conformal field theory (Charlottesville, VA, 2000),
201--225, { Contemp. Math.} {\bf 297}, Amer. Math. Soc.,
Providence, RI, 2002.

\bibitem[42]{M1} A. Milas, Logarithmic intertwining operators and
vertex operators, to appear in {\em Comm. Math. Physics}; {\tt
math.QA/0609306}.

\bibitem[43]{Miy} M. Miyamoto,
Modular invariance of vertex operator algebras satisfying $C\sb
2$-cofiniteness. {\em Duke Math. J.} {\bf 122} (2004), 51--91.

%\bibitem[36]{TK} A. Tsuchiya and Y. Kanie, Fock space representations of
%the Virasoro algebra - Intertwining operators,  { Publ. RIMS} {\bf
%22} (1986), 259-327.

%\bibitem[W1]{W1} W.Wang, Rationality of Virasoro Vertex operator
%algebras, Internat. Math. Res. Notices , Vol {\bf 71}, No.1
%(1993), 197-211.

%\bibitem[37]{W2}  W. Wang,  ${\W}_{ 1+ \infty}$ algebra, ${\W}_3$ algebra, and
% Friedan\--Martinec\--Shenker bosoni\-zation,
%{ Comm. Math. Phys.} {\bf  195} (1998), 95-111.
%

%\bibitem[38]{W3}  W. Wang,
% \ Classification of irreducible modules of ${\W}_3$ algebra
%with  $c=-2$, { Comm. Math. Phys.} {\bf 195} (1998),  113-128.
%

\bibitem[44]{Zhu}
Y.-C. Zhu, Modular invariance of characters of vertex operator
algebras, J. Amer. Math. Soc.  9 (1996), 237-302.

\end{thebibliography}
\end{document}